\documentclass[12pt]{article}
\usepackage{amsmath,amsfonts,amssymb}
\usepackage{mathrsfs}
\usepackage{dsfont}
\usepackage[utf8x]{inputenc}
\usepackage{enumitem}
\usepackage[T1]{fontenc}
\usepackage{lmodern} \usepackage{ucs}
\usepackage{fullpage}

\usepackage{graphicx}
\usepackage{caption}
\usepackage{subcaption}

\usepackage{array}
\usepackage{multicol}
\usepackage{multirow}
\usepackage{color}
\usepackage{authblk}
\usepackage{esint}

\usepackage{tikz}
\usetikzlibrary{decorations}
\usetikzlibrary{decorations.pathreplacing,calc}
\usetikzlibrary{patterns}

\newtheorem{lemma}{Lemma}[section]
\newtheorem{theo}[lemma]{Theorem}
\newtheorem{rmk}[lemma]{Remark}
\newtheorem{proposition}[lemma]{Proposition}

\newtheorem{coro}[lemma]{Corollary}

\newcommand{\hslashslash}{%
  \raisebox{.9ex}{%
    \scalebox{.7}{%
      \rotatebox[origin=c]{18}{$-$}%
    }%
  }%
}

\def\dslash{%
  {%
   \vphantom{d}%
   \ooalign{\kern.05em\smash{\hslashslash}\hidewidth\cr$\mathrm d$\cr}%
   \kern.05em
  }}
  
\newcommand{\QED}{\mbox{}\hfill \raisebox{-0.2pt}{\rule{5.6pt}{6pt}\rule{0pt}{0pt}} \medskip\par}

\newcommand{\R}{\mathbb{R}}

\newcommand{\ds}{\displaystyle}
\newcommand{\ud}{\, {\mathrm{d}}}

\title{Solutions of the 4-species quadratic reaction-diffusion system are bounded and $C^\infty$-smooth, in any space dimension
}

\author[1]{M. Cristina Caputo}

\author[2]{Thierry Goudon\thanks{ {\tt thierry.goudon@inria.fr}}}

\author[1]{Alexis F. Vasseur\thanks{ {\tt vasseur@math.utexas.edu}}}

\affil[1]{\small Department of Mathematics,
University of Texas at Austin

1 University Station C1200
Austin, TX, 78712-0257
USA }

\affil[2]{\small Universit\'e C\^ote d'Azur, Inria,  CNRS, LJAD,

Parc Valrose, F-06108 Nice, France}

\begin{document}
\maketitle

\abstract{We establish the boundedness of solutions of reaction-diffusion systems with qua\-dra\-tic (in fact slightly super-quadratic) 
reaction terms that satisfy a natural entropy dissipation property,  in any space dimension $N>2$. This bound imply the $C^\infty$-regularity of the solutions. 
This result extends the theory which was restricted to the two-dimensional case. The proof heavily uses De Giorgi's iteration scheme, which allows us to 
obtain local estimates. The arguments relies on  duality reasonings  in order to obtain new estimates on the total mass of the system, both in $L^{(N+1)/N}$ norm and in a suitable weak norm. The latter uses $C^\alpha$ regularization properties for parabolic equations.}

\vspace*{.5cm}
{\small
\noindent{\bf Keywords.}
Reaction-diffusion systems. Global regularity. Blow-up methods. \\[.4cm]

\noindent{\bf Math.~Subject Classification.} 
35K45, 
35B65, 
35K57. 
}

\section{Introduction}

This paper is mainly concerned with the following system of reaction-diffusion equations
\begin{equation}\label{eqA}\begin{array}{l}
\partial_t a_i-\nabla\cdot(D_i\nabla a_i)=Q_i(a),\qquad i\in\{1,2,3,4\},\ t\geq 0,\ x\in\mathbb R^N,
\\[.3cm]
Q_i(a)=(-1)^{i+1}(a_2a_4-a_1a_3),
\end{array}\end{equation}
with initial condition
\begin{equation}\label{eqAci}
a\big|_{t=0}=a^0=(a_1^0,a_2^0,a_3^0,a_4^0).
\end{equation}
This system arises in chemistry where four species interact according to the reactions
\[A_1 + A_3 \rightleftharpoons  A_2+A_4,\]
the unknowns $(t,x)\mapsto a_i(t,x)$ in \eqref{eqA} being the local mass concentrations
of the species labelled by  $i\in\{1,2,3,4\}$: $\int_{\mathbb R^N} a_i(t,x)\ud x$ is interpreted as the mass of the constituent $i$ at time $t$.
It is thus physically relevant to consider  initial data $a_i^0$ which are  non negative integrable functions.
The reactants are  subjected to space diffusion and 
the diffusion coefficients
depend on the considered species. 
In full generality, $D_i$ can be a function of the space variable with values in the space of $N\times N$ matrices. Throughout this paper, we restrict to the case of scalar and constant matrices
\[D_i(x)=d_i \mathbb I,\qquad \textrm{$d_i>0$ constant}\]
with coefficients that satisfy
\begin{equation}
\label{hypcoef}
\textrm{
$ 0<\delta_\star\leq d_i \leq \delta^\star.$}
\end{equation}
Assuming that the initial data are smooth, say $a_i^0\in C^\infty(\mathbb R^N)$, 
existence-uniqueness of smooth and non-negative solutions for \eqref{eqA}--\eqref{eqAci} can be justified at least on a small time interval, by using 
a standard fixed point reasoning (see for instance \cite[Proposition~A.2]{GoVa} or \cite[Lemma~1.1]{Pie}).
Global existence of weak solutions is established in \cite{DFPV}.
We address the question of the boundedness
of the solutions, which will imply that solutions are globally defined and remain infinitely smooth  \cite[Proposition~A.1]{GoVa}.

The difficulty comes from the fact we are dealing with  different diffusion coefficients.
As already noticed in \cite{GoVa}, the question becomes trivial when all the $D_i$'s vanish: in this case, we are 
concerned with a mere system of ODE which clearly satisfies
a maximum principle.
The answer is also immediate when all the diffusion coefficients are equal to the same constant
$d_i=\delta_\star$. Indeed, in this situation, the total mass
\[
M(t,x)=\ds\sum_{i=1}^4 a_i(t,x)\]
satisfies the heat equation
$\partial_t M=\delta_\star\Delta M$, which, again, easily leads to a maximum principle.
In the general situation, one may wonder whether or not the system 
has the explosive behavior of non linear heat equations \cite{We}.
Counter--examples of systems with polynomial non linearities presented in \cite{PiSc2} show that this question is relevant and
non trivial, see also \cite[Theorem~4.1]{Pie}.
We refer the reader to the survey \cite{Pie} for a general presentation of the problem, further references, and many deep comments 
on the mathematical difficulties raised by such systems.
\\

Two properties are crucial for the analysis of the problem. First of all, system \eqref{eqA}
conserves mass
\begin{equation}\label{masscons}
\ds\frac{\ud}{\ud t}\ds\sum_{i=1}^4 \ds\int_{\mathbb R^N} a_i\ud x=0.
\end{equation}
Second of all, it dissipates entropy: 
\begin{equation}\label{ent_d}
\ds\sum_{i=1}^4 Q_i(a)\ln (a_i)=-(a_2a_4-a_1a_3)\ln\Big(\ds\frac{a_2a_4}{a_1a_3}\Big)\leq 0.
\end{equation}
These properties  suggest to consider more general systems, involving more reactants and possibly more intricate non linearities.
To be more specific, we extend the discussion to systems that read
\begin{equation}\label{eqA2}\begin{array}{l}
\partial_t a_i-\nabla\cdot(D_i\nabla a_i)=Q_i(a),\qquad i\in\{1,...,p\},\ t\geq 0,\ x\in\mathbb R^N,
\\[.3cm]
Q_i: a\in\mathbb R^p\longmapsto \mathbb R^p,
\end{array}\end{equation}
endowed with the initial condition
\begin{equation}\label{eqA2ci}
a\big|_{t=0}=a^0=(a_1^0,...,a_p^0),
\end{equation}
where the reaction term 
fulfils the following conditions
\begin{enumerate}
\item [h1)] there exists $\mathscr Q>0$ and $q>0$ such that for any $a\in \mathbb R^p$ 
and $i\in \{1,...,p\}$, we have $|\nabla_a Q_i(a)|\leq \mathscr Q |a|^{q-1}$,
\item [h2)] for any $i\in \{1,...,p\}$, if $a_i\leq 0$ then $Q_i(a)\leq 0$,
\item [h3)] $\ds\sum_{i=1}^p Q_i(a)=0,$
\item [h4)]  $\ds\sum_{i=1}^p Q_i(a)\ln(a_i)\leq 0.$
\end{enumerate}
Assumption h1) governs the growth of the non linearity. In what follows, we will be concerned with 
quadratic and super-quadratic growth: $q\geq2$ (but $q$ is not necessarily assumed to be an integer). Assumption h2) relies on the preservation of non negativity of the solutions, and it is thus physically relevant.
Assumptions h3) and h4) imply mass conservation and entropy dissipation, respectively.
Note that the entropy dissipation actually provides an estimate on (nonlinear) derivatives of the unknown since it leads to
\begin{equation}\label{ent_diss}
\ds\frac{\ud}{\ud t}\ds\sum_{i=1}^p \ds\int_{\mathbb R^N} a_i\ln(a_i)\ud x
+4\delta_\star\ds\sum_{i=1}^p \ds\int_{\mathbb R^N} |\nabla\sqrt {a_i}|^2\ud x\leq 0.
\end{equation}
In view of h3) and h4), it is thus natural to consider initial data such that
\begin{equation}\label{hyp_ci}\begin{array}{l}
a_i^0:x\in \mathbb R^N\longmapsto a_i^0(x)\geq 0,
\\[.3cm]
\ds\sup_{i\in\{1,...,p\}} \ds\int_{\mathbb R^N} a_i^0(1+\ln (a_i^0)+ |x|)\ud x=\mathscr M^0<\infty. 
\end{array}
\end{equation}
We refer the reader to Proposition~\ref{P1} below for a more precise statement in terms of a priori estimate.
It means that the initial concentrations have finite mass and entropy.
The moment condition 
controls the spreading of the mass. However, while \eqref{ent_diss} has a clear physical meaning, 
it does not provide enough estimates for the analysis of the problem:
note that with $u, u\ln (u) \in L^1$ and $\nabla\sqrt u\in L^2$, it is still not clear how the nonlinear term  $Q(u)$ can make sense in $\mathscr D'$~!
For this reason, a notion of  \emph{renormalized} solutions is introduced in \cite{Fis}, and existence of solutions in this framework can be established.
\\

In the specific quadratic and two-dimensional case ($q=2$, $N=2$) the question is fully answered  in
\cite{GoVa}: starting from $L^\infty\cap C^\infty$ initial data, the solution remains bounded and smooth 
and the problem is globally well-posed. In fact \cite{GoVa} proves a \emph{regularizing} effect:
with data satisfying \eqref{hyp_ci} only, the solution becomes \emph{instantaneously} bounded and smooth, 
which implies global well-posedness.  
The proof in \cite{GoVa} relies on De Giorgi's approach \cite{DeG}; it uses entropy dissipation, see \eqref{ent_diss},
to get a non linear control on level sets of the solution, which eventually leads to the $L^\infty$ bound.
The result is extended for  higher space dimensions
in \cite{CDF} which handles the quadratic case when the diffusion coefficients are 
close enough to the same  constant
(how small the distance between the   $d_j$'s should be depends on the space dimension, in a explicit way),
and 
in \cite{CaVa}, which handles subquadratic non linearities ($q<2$ in h1), non necessarily integer).
Two ingredients
are crucial in the approach of \cite{CaVa}:
\begin{itemize}
\item First, \cite{CaVa} uses systematically rescaled quantities
\begin{equation}\label{scal}
a^{(\epsilon)}_i(s,y)=\epsilon^{2/(q-1)}\ a_i(t+\epsilon^2 s,x+\epsilon y)
\end{equation}
with $\epsilon>0$: $a^{(\epsilon)}$ satisfies the same evolution equation as $a$. 
Note that in the quadratic case ($q=2$), for $N=2$, the rescaling leaves invariant the natural norms of the problem
$\|a\|_{L^\infty(0,\infty;L^1(\mathbb R^2))}$ and $\|\nabla\sqrt a\|_{L^2((0,\infty)\times\mathbb R^2)}$.
\item Second, the parabolic regularity is obtained by adapting De Giorgi's techniques, and by working  with a certain norm of the rescaled unknown
which becomes small as $\epsilon\to 0$.
It turns out that the necessary estimate holds in a weak sense. Namely, one has to consider the set of distributions
\[\textrm{$T\in \mathscr D'((0,T)\times\mathbb R^N)$ such that $T=\Delta \Phi$, 
 with $\Phi\in L^\infty((0,\infty)\times\mathbb R^N)$}.\]
 The corresponding rescaled norm behaves like $\mathscr O(\epsilon^{(4-2q)/(q-1)})$, which indeed tends to 0 as $\epsilon \to 0$ for 
 subquadratic non linearities $q<2$. The idea of using such a weak norm also appeared in  the regularity analysis for the 
 Navier-Stokes equation \cite{Va}.
 We also refer the reader to \cite{CafVa, Va0},  for further applications of De Giorgi's techniques to the analysis of fluid mechanics systems
 and to \cite{AAG, GU} for the study of models for populations dynamics governed by ``chemotactic-like'' mechanisms.
 This approach is also useful for the analysis 
  of the preservation of bounds by numerical schemes when solving non linear convection-diffusion systems \cite{CMV}.
  In the reasoning adopted in \cite{CaVa}, a special role is played by the total mass
 $M=\sum_{i=1}^p a_i$ which satisfies the diffusion equation
 \begin{equation}\label{eqMass}
 \begin{array}{l}
 \partial_t M - \Delta (dM)=0,
 \\[.3cm]
 d(t,x)=\ds\frac{\ds\sum_{i=1}^p d_i a_i(t,x)}{\ds\sum_{i=1}^p a_i(t,x)},
 \end{array}\end{equation} 
 where, by virtue of \eqref{hypcoef}, the diffusion coefficient $d$ satisfies
 \[0<\delta_\star\leq d(t,x)\leq \delta^\star.\]
\end{itemize}
This relation can be used to establish, through an elegant duality argument, an estimate in $L^2((0,T)\times\mathbb R^N)$, see \cite{PiSc2} and \cite{DFPV}.
This estimate is a key for proving the global existence of weak solutions for the quadratic problem \eqref{eqA}--\eqref{eqAci} in \cite{DFPV}: at least, it is  worth  pointing out that with this $L^2$ estimate
the right hand side $Q_i(a)$ in \eqref{eqA} makes sense, while the estimates based on the mass conservation and entropy dissipation were not enough.
However, the $L^2$ estimate does not schrink the rescaled solutions $a^{(\epsilon)}$ as $\epsilon\to 0$ and it is thus not enough to provide global boundedness and regularity.
This is where we can take advantage of using  a weak norm.

In the present work, we wish to fill the gap in the boundedness theory and to provide a complete answer for the quadratic case in \emph{any} dimension.
In fact, our analysis also covers higher non linearities, but with a non explicit condition on the growth condition.
Our main results 
state as follows.

\begin{theo}\label{theo_main} 
Let $N\in\mathbb N$, with $N\geq 3$.
For any  initial data $a^0=(a^0_1,a^0_2,a^0_3,a^0_4)$ in $\big(C^\infty(\R^N)\cap L^\infty(\R^N)\big)^4$ 
such that $a_i(x)\geq0$ for any $x\in \R^n$ and $i\in\{1,...,4\}$, 
 there exists a unique,  globally defined,
 solution $a= (a_1,a_2,a_3,a_4)$ to \eqref{eqA}--\eqref{eqAci} which is non negative, 
 bounded on $[0,T]\times\mathbb R^N$ for any $0<T<\infty$, and
 $C^\infty$-smooth. \end{theo}

\begin{theo}\label{Theo_principal} 
Let $N\in\mathbb N$, with $N\geq 3$. Consider a system \eqref{eqA2} verifying h1)-h4).
There exists $\nu_0>0$ depending on $N$, $\delta_\star$ and $\delta^\star$ such that if h1) holds with $2\leq q\leq 2+\nu_0\leq 2\frac{N+1}{N}$, then 
for any  non negative $a^0 \in C^{\infty}(\R^N;\R^p)\cap L^{\infty}(\R^N;\R^p)$, 
there
exists a unique,  globally defined,  solution $a$ to \eqref{eqA2}--\eqref{eqA2ci} which is non negative, bounded on $[0,T]\times\mathbb R^N$ for any $0<T<\infty$,
and $C^\infty$-smooth.
\end{theo}

 Theorem \ref{theo_main} thus appears as a consequence of  Theorem \ref{Theo_principal}. 
The extra power $\nu_0$ allowed on the nonlinearities depends on $N$, $\delta_\star$ and $\delta^\star$ in a non explicit way and our method does not provide any precise estimate. It seems unlikely that it can correspond to a physically relevant threshold. The problem of regularity remains open for higher nonlinearities. 
The proof still follows the De Giorgi strategy, and relies on a refinement of the weak norm estimate obtained in \cite{CaVa}
(which, though, remains a crucial ingredient of the proof). To be more specific,
we are going to upgrade the $L^\infty$ estimate to a $C^\alpha$ estimate, working with the set of distributions
\[\textrm{$T\in \mathscr D'((0,T)\times\mathbb R^N)$ such that $T=\Delta \Phi$, 
 with $\Phi\in L^\infty(0,\infty;C^\alpha(\mathbb R^N))$}\]
for a certain regularity coefficient $0<\alpha\leq 1$.
This is combined with a $L^{(N+1)/N}$ estimate on the total mass, obtained through a duality argument.
This argument  is directly inspired by the derivation of elliptic estimates by Fabes and Stroock \cite{FaSt} and it appears as a dual version of the 
 Alexandrof-Bakelman-Pucci-Krylov-Tso  (ABPKT) estimate
\cite{Ale,Bak,Puc,Kry,Tso}.
We point out that, contrarily to the approach in \cite{CaVa}, 
we do not use here the bounds derived from the entropy dissipation \eqref{ent_diss}.

The paper is organized as follows. In Section \ref{S:main}, we give an overview of the main steps of the proof.
Section \ref{S:3} is concerned with the weak estimate on the total mass. It relies on a H\"olderian regularity analysis for parabolic equations. This is combined with 
a duality argument which uses crucially the non negativity of the solution.
Section \ref{S:4} is devoted to a complementary 
estimate in a suitable Lebesgue space, which, again, relies on a duality approach. Section \ref{S:fin} explains how the arguments combine to end the proof of the main results.

\section{Main steps of the proof}
\label{S:main}

\subsection{A priori estimates; boundedness, global existence and regularity of the solutions}

In what follows, we are going to establish several a priori estimates  satisfied by the solutions
of \eqref{eqA2}.
To this end, we will perform  various manipulations such as integrations by parts,
permutations of integrals and derivation, etc.
These manipulations apply to the  smooth solutions of the problem
that can be shown to exist on a small enough time interval, see \cite[Proposition~A.2]{GoVa}.
They equally  apply to solutions of suitable approximations
of the problem \eqref{eqA2}. 
The construction of such an approximation --- by regularizing data, coefficients, 
cutting-off the non linearirities... --- can be a delicate issue in order to 
preserve the structural features of the original equation, 
and to admit a globally defined smooth solution.
We refer the reader on this issue to \cite{DFPV}.
As it will be clear in the forthcoming discussion, the estimates we are going to derive do not depend on the regularization parameter, but only 
on  $N$, $\delta_\star$, $\delta^\star$, and $\mathscr Q$,  $p$, $q$ (see h1)), 
which, eventually, allows us 
to conclude by getting rid of the regularization parameter.
The very first estimate is a direct consequence of the mass conservation and entropy dissipation properties of the system.
The following claim, see \cite[Proposition~2.1]{GoVa}, applies without any restriction on 
the number of species
$p$, the degree of non linearity $q$ nor on the space dimension $N$.


\begin{proposition}[\cite{GoVa}]\label{P1}
Assume h1)-h4).
Let $a_0=(a_1^0,...,a_p^0)$, with non negative components,  satisfy \eqref{hyp_ci}. Then, for any $0<T<\infty$, there exists $0<C(T)<\infty$ such that
\[
\begin{array}{l}
\ds\sup_{0\leq t\leq T}\Big\{\ds\sum_{i=1}^{p}\ds\int_{\mathbb R^N}
a_i\big(1+|x|+|\ln(a_i)|\big)(t,x) \ud x\Big\}
\\
\qquad\qquad
+
\ds\sum_{i=1}^{p}\ds\int_0^T\ds\int_{\mathbb R^N} \big|\nabla
\sqrt{a_i}\big|^2(s,x) \ud x\ud s
+
\ds\sum_{i=1}^{p}\ds\int_0^T\ds\int_{\mathbb R^N}
Q_i(a)\ln (a_i)
\ud x\ud s
\leq C(T).\end{array}
\]
\end{proposition}

The entropy dissipation \eqref{ent_diss} tells us that 
$\sum_{i=1}^{p}\int_{\mathbb R^N}
a_i\ln(a_i)(t,x) \ud x$
is a non increasing function of the time variable. However, this quantity has no sign.
To make this information a useful estimate, involving the non negative quantities
$a_i|\ln(a_i)|
$ we need a control  on the first order space moments $\int_{\mathbb R^N}
|x| a_i(t,x) \ud x
$. We refer the reader to \cite{GoVa} for details.
 This estimate  will not be used in our reasoning; nevertheless the entropy dissipation
 still has a crucial role in the proof of Theorems~\ref{theo_main} and~\ref{Theo_principal}.
 By the way note that the counter examples of systems that produce blow up 
 in \cite{PiSc2} very likely do not satisfy the entropy dissipation property.
\\

As said above, for data in $C^\infty\cap L^\infty(\mathbb R^N)$, we can construct a $C^\infty$ and bounded solution defined on a small enough interval.
Let $T_{\mathrm{max}}$ be the lifespan of such a solution. 
Standard bootstrapping arguments tell us that if $T_{\mathrm{max}}<\infty$ then we have
 $$
 \limsup_{t\to T_{\mathrm{max}}}\| a (t,\cdot)\|_{L^\infty(\R^N)}=+\infty.
 $$
In what follows, we are going to  obtain a uniform bound satisfied by  $\|a(t,\cdot)\|_{L^\infty(\R^N)}$ on the time interval $[0,T_{\mathrm{max}})$,
depending  only on $T_{\mathrm{max}}$ and the assumptions on the data,
which thus  contradicts the occurrence of a blow-up of the solution in finite time. 
Therefore, the $L^\infty$ estimate implies that the lifespan of the solutions of \eqref{eqA2}--\eqref{eqA2ci} is infinite.
Moreover, boundedness also implies the regularity of the solution, by a bootstrap argument, see \cite[Proposition~A.1]{GoVa}.

\subsection{The key intermediate statements}

The main ingredient consists in showing that
the local boundedness can be obtained from 
 a local estimate in $L^r$, with $r>1$, see  \cite[Proposition~4]{CaVa}. 
 We thus work on balls 
 \[B_\rho=\big\{x\in\mathbb R^N,\ |x|\leq \rho\big\}
.\]

\begin{lemma}[De Giorgi type Lemma, \cite{CaVa}]\label{lemm_de giorgi}
We suppose that  $2\leq q< \frac{2(N+1)}{N}$. We also suppose that h1)-h4) holds.
Let $a $ be a non negative solution to \eqref{eqA2} on $(-1,0)\times {B}_1$. Then,
for any $r>1$, there exists a universal constant $\delta_r>0$ such that, if  $a=(a_1,\cdot\cdot\cdot, a_p)$
verifies
$$ \ds\sum_{i=1}^p \|a_i\|_{L^{ r}((-1,0)\times B_1)}\leq \delta_r,$$
then, $0\leq a_i(0,0)\leq 1$, for $i\in\{1,...,p\}$.
\end{lemma}

The proof relies on De Giorgi's techniques \cite{DeG} (see also \cite{Alikakos} for a related approach). For the sake of completeness we describe the main steps in Appendix~\ref{a:L3};
it is also important to detail this proof since this is where the entropy dissipation plays a central role.
At first sight this information does not look very useful since the natural estimates 
for \eqref{eqA2}--\eqref{eqA2ci} in Proposition~\ref{P1} do not involve $L^r$ norms for an exponent $r$ larger than 1.
However, we will be able to identify further estimates, 
that shrink 
for the rescaled solutions \eqref{scal} as $\epsilon\to 0$.
Thus, for $\epsilon$ small enough the rescaled solution fulfils 
the criterion in Lemma~\ref{lemm_de giorgi}.

\begin{lemma}\label{lemm_de giorgi2} 
There exists $\epsilon_0>0$ and $\nu_0>0$ depending on $N$, $\delta_\star$ and $\delta^\star$ such that if h1) holds with $2\leq q \leq 2+\nu_0\leq 2\frac{N+1}{N}$,
then  for all $0<\epsilon\leq\epsilon_0$ we have
$$\sum_{i=1}^{p} \|{ a}^{(\epsilon)}_i\|_{L^{(N+1)/N}((-1,0)\times B_1)}\leq \delta$$
with $\delta=\delta_{(N+1)/N}$ as defined  in Lemma~\ref{lemm_de giorgi}.
\end{lemma}

\noindent
Coming back to the original variables, we obtain the $L^\infty$ estimate.

\begin{coro}\label{unibound} Let $\epsilon_0$ be defined in  Lemma~\ref{lemm_de giorgi2}. Then, for all $\frac{T_{\mathrm{max}}}{2}<t<T_{\mathrm{max}}$, we have
$$\sum_{i=1}^{p} \|a_i(t,\cdot)\|_{L^{\infty}(\R^N)}\leq \epsilon_0^{-2/(q-1)}.$$
\end{coro}

\noindent
{\bf Proof.}
Pick $x_0$  in $\R^N$ and  $t_0\in (\frac{T_{\mathrm{max}}}{2},T_{\mathrm{max}})$. 
Applying Lemma~\ref{lemm_de giorgi} to $a^{(\epsilon_0)}$ yields
$$0\leq \sum_{i=1}^{p}  a_i(t_0,x_0)=\epsilon_0^{\frac{-2}{q-1}}\sum_{i=1}^{p} a_i^{(\epsilon_0)}(0,0)\leq \epsilon_0^{\frac{-2}{q-1}}.$$
\QED

Having this statement at hand allows us to conclude the proof of Theorem~\ref{Theo_principal}.
Let $2\leq q\leq 2+\nu_0\leq 2\frac{N+1}{N}$. Let $a=(a_1,..., a_p)$ be a solution to~\eqref{eqA2}--\eqref{eqA2ci}, and let $T_{\mathrm{max}}$ be the lifespan of $a$. Assume that $T_{\mathrm{max}}$ is finite. Then, for  each $i\in \{1,...,p\}$,   Corollary~\ref{unibound} tells us that $a_i(t,\cdot)$ is  uniformly bounded for all $\frac{T_{\mathrm{max}}}{2}<t<T_{\mathrm{max}}$
and thus the sup norm does not blow up as $t\to T_{\mathrm{max}}$. This contradicts  the fact that $T_{\mathrm{max}}$ is the maximal time of existence of a smooth solution of \eqref{eqA2}--\eqref{eqA2ci}.
\QED

Therefore the cornerstone of the proof consists in proving Lemma~\ref{lemm_de giorgi2} and identifying the specific role payed by the norm $L^{(N+1)/N}$.
The argument is two-fold and it uses the diffusion equation \eqref{eqMass}
satisfied by the total mass
$M(t,x)=\sum_{i=1}^p a_i(t,x)$.
On the one hand, we shall show that the  norm $L^{(N+1)/N}$ of $M$ can be controlled by means of the norm  $L^\infty(0,\infty;L^1(\mathbb R^N))$.
On the other hand, we shall obtain a new estimate on a \emph{weak norm} of $M$, 
which will allow us to conclude that 
\[\ds\lim_{\epsilon\to 0} \|M^{(\epsilon)}\|_{L^\infty(0,\infty;L^1(\mathbb R^N))}  =0,\qquad
\textrm{with $M^{(\epsilon)}(s,y)=
\epsilon^{2/(q-1)}\ M(t+\epsilon^2 s,x+\epsilon y)$}.\]
This analysis is based on duality arguments and regularization properties of parabolic equations.
Accordingly, we can conclude to the shriking as $\epsilon\to 0$ of the $ L^{(N+1)/N}$ norm of the rescaled solutions.

\subsection{Preliminary comments}

The  De Giorgi approach leads us to construct sequences, based on energy-entropy estimates, where the parameter of the sequence 
controls level sets of the solution and space-time localization. 
Roughly speaking, we obtain a non linear control of the $k$th level by the $(k-1)$th level.
We can finally conclude to a local property of the solution by using the following simple result.

\begin{lemma}\label{2bete}
Let $\big( u_n\big)_{n\in\mathbb N}$ be a sequence of non negative real numbers.
We suppose that it satisfies, for any $n\in \mathbb N \setminus\{0\}$,
$$u_n\leq \Lambda ^n u_{n-1}^{\gamma}$$
where $\Lambda, \gamma>1$.
Then, there exists $\kappa>0$ such that, if $0\leq u_0\leq \kappa$, then
$\lim_{n\to \infty} u_n=0$.
\end{lemma}

\noindent
{\bf Proof.}
We set $v_n=\ln (u_n)$ which satisfies
\[v_n\leq n\ln(\Lambda)+ 
\gamma  v_{n-1},\]
and thus
\[
v_n\leq \ln(\Lambda)\ds\sum_{j=0}^n j\gamma^{n-j} 
+v_0 \gamma^n
\leq \gamma^n\ln (\Lambda^{F(\gamma)} u_0)
\]
with 
$$F(\gamma)=\ds\frac1\gamma\ds\sum_{j=0}^\infty  j\Big(\ds\frac1\gamma\Big)^{j-1}
=\ds\frac1\gamma \ \ds\frac{\ud}{\ud x}\Big(\ds\frac{1}{1-x}\Big)\Big|_{x=1/\gamma}=
\ds\frac1\gamma\Big(\ds\frac{1}{1-1/\gamma}\Big)^2.  $$
Therefore  $v_n$ tends to $-\infty$, and $u_n$ tends to 0, as $n\to \infty$ provided $u_0$ is small enough.
\QED

\section{Weak norm estimates on the total mass and shrinking of the rescaled total mass}
\label{S:3}


Our approach relies on the following statement.

\begin{proposition}\label{p:1}
Let $\Phi:(0,T)\times \mathbb R^N\rightarrow \mathbb R$ such that
\begin{itemize}
\item[a)] $\Phi$ lies in $L^\infty((0,T)\times\mathbb R^N)$; 
\item[b)] $\Delta \Phi=M\geq 0$;
\item[c)]  $\Phi$ satisfies $\partial_t \Phi-d\Delta \Phi=0$ on  $(0,T)\times \mathbb R^N$, with 
a coefficient  $d:(0,T)\times \mathbb R^N\rightarrow \mathbb R$ verifying 
$0< \delta_\star \leq d(t,x)\leq \delta^\star<\infty$   for a.\,e.\,$(t,x)\in (0,T)\times \mathbb R^N$.
\end{itemize} 
Then, there exists $\alpha\in (0,1]$ such that 
$\Phi\in C^{[\alpha/2,\alpha]}([t_0,T]\times \mathbb R^N)$ for any $t_0>0$, which means that
we can find $C>0$ such that, 
for any $(t,x)\in [t_0,T]\times \mathbb R^N$ and  $(\tau,h)\in \mathbb R\times \mathbb R^N$
with $t+\tau \geq t_0$, we have
\[
\ds\frac{|\Phi(t+\tau,x+h)-\Phi(t,x)|}{|\tau|^{\alpha/2}+|h|^\alpha}\leq C \|\Phi\|_{L^\infty}.\]
\end{proposition}

\noindent
This H\"older regularity estimate for non conservative parabolic equations dates back to Krylov-Safonov \cite{SK1,SK2}.
In fact, the result of  \cite{SK1,SK2} does not need the sign property b).
However, as it will be explained below, 
this sign property naturally appears for the system under consideration, and it plays 
a further crucial role throughout the analysis.
Let us explain the interest of this statement for our purpose.
As said above the total mass $M$ satisfies the diffusion equation \eqref{eqMass}.
Of course, by definition, $M$ is a non negative function which lies in $L^\infty(0,\infty;L^1(\mathbb R^N))$.
Let $\Phi$ satisfy $\Delta \Phi =M\geq 0$.
Since  $d(t,x)$ is bounded above by $\delta^\star$, $\Phi$ also satisfies the evolution equation
\[
\partial_t \Phi-\delta^\star \Delta \Phi=(d-\delta^\star)\Delta \Phi=(d-\delta^\star)M\leq 0.\]
This observation is the cornerstone of the analysis performed in \cite{CaVa}.
In particular, we will make use of the following crucial property established in \cite[Proposition~11 \& Corollary~12]{CaVa}.
\begin{proposition}\label{P:CaVa}
Let $N\in\mathbb N$, with $N\geq 3$.
Let $\Phi=\Delta^{-1}M$ with $M$ the total mass associated to a solution of \eqref{eqA2}.
Then, we have
\[
\|\Phi\|_{L^\infty((0,T)\times\mathbb R^N)}\leq \|\Phi(0,\cdot)\|_{L^\infty(\mathbb R^N)}\leq K_N\
\| M(0,\cdot)\|_{L^\infty(\mathbb R^N)}^{1-2/N}\ \| M(0,\cdot)\|_{L^1(\mathbb R^N)}^{2/N}
,\]
where $K_N>0$ is a certain universal constant, which only depends on the space dimension.
\end{proposition}

Proposition~\ref{p:1} thus strengthens \cite{CaVa}'s results in the sense that it provides, beyond the $L^\infty$ 
estimate on $\Phi$, a H\"older-regularity estimate. Since the estimate in Proposition~\ref{P:CaVa} is not evident at first sight, we give the main steps of the proof in Appendix~\ref{app:CaVa} for the sake of completeness.
We shall use the following consequence of Proposition~\ref{p:1}, which is precisely the estimate that allows us to go 
beyond the subquadratic non linearities dealt with in \cite{CaVa}.

\begin{lemma}\label{l:L1shr}
Let $M$ be a non negative solution of \eqref{eqMass}, and let $\Phi=\Delta^{-1}M$. Let $t\geq t_0>0$ and $x\in\mathbb R^N$.
For $\epsilon>0$, we set
$M^{(\epsilon)}(s,y)=\epsilon^{2/(q-1)}M(t+\epsilon^2 s,x+\epsilon y)$.
We suppose that $M^{(\epsilon)}$ lies in $L^\infty(-4,0;L^1(\mathbb R^N))$.
Then, there exists $c>0$ and $0<\alpha\leq 1$, depending only on $N$, $\delta_\star$ and $\delta^\star$, such that
for any $0<\epsilon\leq \sqrt {t_0}/2$,
\[
\ds\sup_{-4\leq s\leq 0}\ds\int_{B_2} M^{(\epsilon)}(s,y)\ud y\leq c\  \|\Phi\|_{L^\infty}\ \epsilon^{\alpha-2+2/(q-1)}\ 
.\]
\end{lemma}

\noindent
{\bf Proof.}
Let $\zeta\in C^\infty_c(\mathbb R^N)$ be such that $\mathrm{supp}(\zeta)\subset B_2$ and $\zeta(x)=1$ for any $x\in B_1$.
Since $M^{(\epsilon)}\geq 0$, we get
\[\begin{array}{lll}
\ds\int_{B_1} M^{(\epsilon)}(s,y)\ud y&\leq&
\ds\int_{B_2} \zeta M^{(\epsilon)}(s,y)\ud y=
\ds\int_{B_2} \zeta \Delta \Phi^{(\epsilon)}(s,y)\ud y
\\
[.3cm]
&\leq &\ds\int_{B_2}\Delta \zeta(s,y) \big( \Phi^{(\epsilon)}(s,y)-\Phi^{(\epsilon)}(0,0)\big) \ud y.
\end{array}\]
By virtue of Proposition~\ref{p:1}, we can write
\[\begin{array}{lll}\ds\int_{B_1} M^{(\epsilon)}(s,y)\ud y
&\leq& 
\epsilon^{-2+2/(q-1)}\ds\int_{B_2}\Delta \zeta(y) \big( \Phi(t+\epsilon^2 s,x+\epsilon y)-\Phi(t,x)\big) \ud y
\\[.3cm]
&\leq& C\|\zeta\|_{W^{2,\infty}(\mathbb R^N)}  \|\Phi\|_{L^\infty} \epsilon^{\alpha-2+2/(q-1)} 
\end{array}\]
for any $s\in (-4,0)$ and $0<\epsilon^2<t_0/4$.
\QED

As indicated above the H\"older estimate 
in Proposition~\ref{p:1} is due to \cite{SK1,SK2}.
For the sake of completeness, we provide here an alternative proof, which, however, 
uses the additional assumption b). The interest of this proof is that it entirely relies on energy estimates and De Giorgi's methods.
Since the result stated  in Proposition~\ref{p:1} is standard, the remaining of this Section can be safely skipped
 by the reader not interested in such an alternative proof
 (the original proof relies on a probabilistic interpretation of the equation and uses arguments from the  theory of diffusion 
processes).\\

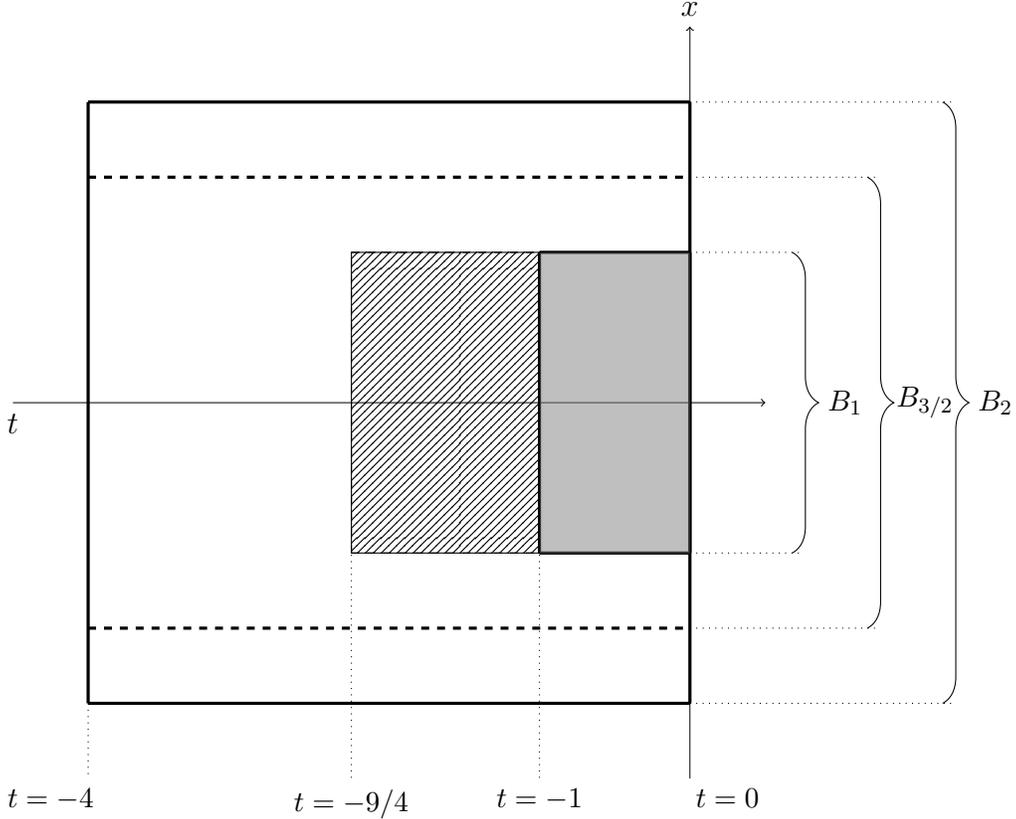
\begin{figure}[!ht]
\begin{center}
\begin{tikzpicture}
\draw [very thin,->]  (0,-5) -- (0,5);
\draw (0,5) node[above] {$x$} ;
\draw [very thin,->]  (-9,0) -- (1,0);
  \draw (-9,0) node[below] {$t$};

   \draw [decorate,decoration={brace,amplitude=10pt},xshift=-4pt,yshift=0pt]
(1.5,2) -- (1.5,-2)node [black,midway,xshift=20pt] {\small $B_1$};

 \draw [decorate,decoration={brace,amplitude=10pt},xshift=-4pt,yshift=0pt]
(3.5,4) -- (3.5,-4)node [black,midway,xshift=20pt] {\small $B_2$};

 \draw [decorate,decoration={brace,amplitude=10pt},xshift=-4pt,yshift=0pt]
(2.5,3) -- (2.5,-3)node [black,midway,xshift=20pt] {\small $\ B_{3/2}$};
  
  \draw [very thick]  (-8,-4) -- (-8,4) ;
  \draw [very thick]  (-8,4) -- (0,4) ;
   \draw [very thick]  (0,4) -- (0,-4) ;
   \draw [very thick]  (-8,-4) -- (0,-4) ;
  
  \draw [dotted] (-8,-4) -- (-8,-5);
   \draw [dotted] (-4.5,0) -- (-4.5,-5);
    \draw [dotted] (-2,0) -- (-2,-5);
     \draw [dotted] (0,-3) -- (2.5,-3);
    \draw [dotted] (0,3) -- (2.5,3);
    
    \draw [dashed, very thick] (-8,-3) -- (0,-3);
    \draw [dashed, very thick] (-8,3) -- (0,3);
    
    \draw [dotted] (0,-4) -- (3.5,-4);
   \draw [dotted] (0,4) -- (3.5,4);
    \draw [dotted] (0,-2) -- (1.5,-2);
     \draw [dotted] (0,2) -- (1.5,2);
  
   \draw (-8.5,-5) node[below] {{\small  {$t=-4$}}};
 \draw (-4.5,-5) node[below] {{\small  {$t=-9/4$}}};
 \draw (-2,-5) node[below] {{\small  {$t=-1$}}};
  \draw (0.5,-5) node[below] {{\small  {$t=0$}}};
  
  \draw [very thick]  (-2,-2) -- (-2,2) ;
  \draw [very thick]  (-2,2) -- (0,2) ;
   \draw [very thick]  (0,2) -- (0,-2) ;
   \draw [very thick]  (-2,-2) -- (0,-2) ;

   \draw [dashed]  (-4.5,-2) -- (-4.5,2) ;
  \draw [dashed]  (-4.5,2) -- (0,2) ;
   \draw [dashed]  (0,2) -- (0,-2) ;
   \draw [dashed]  (-4.5,-2) -- (0,-2) ;

  
\coordinate[] (A) at (-2,-2);
\coordinate[] (B) at (0,-2);
\coordinate[] (C) at (0,2);
\coordinate[] (D) at (-2,2);
\coordinate[] (E) at (-4.5,2);
\coordinate[] (F) at (-4.5,-2);
%
%
%
\fill[gray, opacity=0.5, draw=black] 
(A) -- (B) -- (C) -- (D) -- cycle;

\fill [pattern=north east lines,draw=black]
(D) -- (E) -- (F) -- (A) -- cycle;

%
%
      
 \end{tikzpicture}
\caption{
The domains $Q_2$ (the largest box), $\widetilde Q$ (the dashed box) and $Q_1$ (the grey box)}
\label{Fig:TimeShift}\end{center}
\end{figure}

Here and below, given $\rho>0$, with $B_\rho$ the ball $\{x\in\mathbb R^N, |x|\leq \rho\}$, we denote
\[Q_\rho=(-\rho^2,0)\times B_\rho.\]
In fact, we shall work within $Q_2$, considered as a reference domain. From an equation satisfied on $Q_2$ we wish to establish qualitative properies on a smaller domain, say $Q_1$ or $Q_{1/2}$.
It is also convenient to introduce the domain 
\[\widetilde Q=(-9/4,-1)\times B_1.\]
We refer the reader to Fig.~\ref{Fig:TimeShift}; having the  picture of the subdomains of $Q_2$
 might be helpful in following the arguments.
\\

The argument for proving Proposition \ref{p:1} relies on a technical lemma that controls oscillations.
From now on, for a function $\varphi$ defined on $\Omega\subset \mathbb R^d$, we set
\[
\mathrm{osc}(\varphi,\ \Omega)=\ds\sup_{x\in \Omega} \varphi(x)-\ds\inf_{x\in \Omega} \varphi(x)
.\]

\begin{lemma}[Decay of oscillations]\label{l:osc}
Let $\Phi$ satisfy the assumptions of Proposition~\ref{p:1}.
There exists $\lambda\in (0,1)$, which depends only on $N$ and $\delta_\star$, such that 
\[
\mathrm{osc}\big(\Phi,\ Q_{1/2}\big)\leq \lambda\  \mathrm{osc}\big(\Phi, \ Q_2\big).\]
\end{lemma}

Let us assume temporarily that Lemma \ref{l:osc} holds true.
We pick $(t,x)\in (t_0,T)\times \mathbb R^N$, where $0<t_0<T<\infty$, and we set
\[\Phi_k(t+2^{-2k}s,x+2^{-k}y).\]
where $k\in \mathbb N$ is large enough so that the time variable remains larger than $t_0$ when $-4\leq s\leq 0$; namely, we have 
$k\geq k_0=\ln\big(\frac{t-t_0}{4}\big)\frac{1}{2\ln(1/2)}$.
The function $\Phi_k$
is defined on $Q_2$ and it 
satisfies
\[
\partial_s \Phi_k = d_k \Delta_y \Phi_k\]
where 
\[ d_k(s,y)=d(t+2^{-2k}s,x+2^{-k}y).\]
Moreover, we still have $-1\leq\Phi_k(s,y)\leq +1$.
Applying  Lemma \ref{l:osc}  yields
\[\mathrm{osc}\big(\Phi_k, \ Q_{1/2}\big)\leq \lambda\ \mathrm{osc}\big(\Phi_k, \ Q_2\big)
\]
which rephrases as 
\[\mathrm{osc}\big(\Phi(t+\cdot,x+\cdot), \ Q_{2^{-k-1}}\big)\leq \lambda\  \mathrm{osc}\big(\Phi(t+\cdot,x+\cdot), \ Q_{2^{-k+1}}\big).\]
We deduce that
\[\mathrm{osc}\big(\Phi(t+\cdot,x+\cdot), \ Q_{2^{-k}}\big)\leq \sqrt \lambda^k\times C_0,\qquad
C_0=\ds\frac{2}{\sqrt\lambda^{k_0}}\|\Phi\|_{L^\infty}.\]
(We should bear in mind the fact that $C_0$ depends on $t_0$ through 
the definition of $k_0$ and it is proportional to $\|\Phi\|_{L^\infty}$.)
Let 
$x'\in \mathbb R^N$ and $t'>t_0$; there exists a unique $k\in\mathbb N$ such that 
$x'-x\in B_{2^{-k+1}}\setminus B_{2^{-k}} $, $2^{-2k}\leq| t'-t|\leq 2^{-2(k-1)}$.
It follows that 
\[
\ds\frac{|\Phi(t',x')-\Phi(t,x)|}{|t'-t|^{\alpha/2}+|x'-x|^\alpha}
\leq \ds\frac{C_0}{\sqrt\lambda}\ \big(\sqrt\lambda 2^\alpha)^k.\]
If $0<\sqrt\lambda \leq 1/2$, the right hand side  remains obviously bounded, uniformly with respect to $k$, for any $0<\alpha\leq 1$; otherwise
we  choose 
\[0<\alpha=\ds\frac{\ln(1/\sqrt\lambda)}{\ln(2)}<1.\]
Hence Proposition ~\ref{p:1} follows from Lemma~\ref{l:osc}.
\QED

\noindent
We are thus left with the task of proving  Lemma~\ref{l:osc}.
To this end, we shall apply the following statement.

\begin{proposition}\label{osc}
Let $(t,x)\mapsto v(t,x)$   satisfy 
\begin{itemize}
\item the differential inequality $\partial_t v-\delta^\star\Delta v\leq 0$ on $Q_2$;
\item $-1\leq v(t,x)\leq +1$ on    $Q_2$;
\item $\mathrm{meas}\big(\big\{(t,x)\in  \widetilde Q,\ v(t,x)\leq 0\big\}\big)\geq \mu\ \mathrm{meas}( \widetilde Q)$, for some $\mu>0$
\end{itemize}
Then, there exists $0<\eta<1$ such that 
\[v(t,x)\leq \eta\qquad \textrm{on $Q_{1/2}$}.\]
\end{proposition}

The function 
\[\widetilde \Phi(t,x)=\ds\frac{2}{\mathrm{osc}(\Phi, Q_2)}\Big(
\Phi(t,x)-\ds\frac{\sup_{Q_2}\Phi+\inf_{Q_2}\Phi}{2}\Big)
\]
satisfies the first two assumptions of Proposition \ref{osc}.
Suppose that 
\[\mathrm{meas}\big(\big\{(t,x)\in Q_2,\ \widetilde\Phi(t,x)\leq 0\big\}\big)\geq \ds\frac{\mathrm{meas}(Q_2)}{2}.\]
(Otherwise, we shall apply the same reasoning to $-\widetilde\Phi$.)
Proposition~\ref{osc} tells us that $\widetilde \Phi(t,x)\leq \eta$ on $Q_{1/2}$, which yields 
$\mathrm{osc}(\widetilde \Phi,Q_{1/2})\leq 1+\eta$ (since $\inf_{Q_{1/2}}\widetilde \Phi\geq -1$), and thus
\[\mathrm{osc}(\Phi, Q_{1/2})\leq \ds\frac{1+\eta}{2}\  \mathrm{osc}(\Phi, Q_{2}).\]
It justifies Lemma \ref{l:osc}, with $\lambda=\frac{1+\eta}{2}\in (0,1)$.
\QED

\noindent
The proof of Proposition~\ref{osc} relies on  a series of intermediate statements.

\begin{lemma}\label{distrib}
Let $-\infty<a,b<\infty$ and let $\Omega$ be a smooth bounded domain in $\mathbb R^N$.
We denote $Q=(a,b)\times \Omega$.
\begin{enumerate}[label=(\alph*)]
\item Let $u\in L^\infty(a,b;L^2(\Omega))\cap L^2(a,b;H^1(\Omega))$ such that $$\partial_t u-\delta^\star\Delta u+\mu=0$$
holds in $\mathscr D'(Q)$, with $\mu$ a non negative measure on $Q$. 
Let $F:\mathbb R\rightarrow \mathbb R$ be a non decreasing convex function. 
We assume that $F(0)=0$ and $F\in W^{1,\infty}_{\mathrm{loc}}(\mathbb R)$. Then, there exists a non negative measure $\nu$ such that $v=F(u)$ satisfies
  $\partial_t v-\delta^\star\Delta v+\nu=0$
holds in $\mathscr D'(Q)$
\item Let $v\in L^\infty((a,b)\times \Omega)\cap L^2(a,b;H^1(\Omega))$ be a non negative solution of  $\partial_t v-\delta^\star\Delta v+\nu=0$,
 with $\nu$ a non negative measure on $Q$.
 Then, 
 for any  trial function $\varphi\in C^\infty_c(\Omega)$ 
  there exists $C>0$, which depend only on  $\delta_\star$, $\|v\|_{L^\infty}$ and $\varphi$,
 such that, 
 for a.~e.~$a<s<t<b$, the following energy inequality holds
 $$\begin{array}{l}
 \ds\frac12\ds\int_\Omega v^2(t,x)\varphi^2(x)\ud x +\delta_\star\ds\int_s^t\ds\int_\Omega| \nabla(\phi v)|^2(\tau, x)\ud x\ud \tau
\\
[.3cm]
\qquad\qquad\qquad\qquad \leq  \ds\frac12\ds\int_\Omega v^2(s,x)\varphi^2(x)\ud x + C(t-s).
 \end{array}$$ 
\end{enumerate}
\end{lemma}

\noindent
{\bf Proof.}
Note that $v=F(u)$ also lies in $L^\infty(a,b;L^2(\Omega))\cap L^2(a,b;H^1(\Omega))$, see e.~g.~\cite[Prop.~IX.5]{Brez}.
Item a) follows from the following computation
$$ \partial_t F(u)= -F'(u)\mu + F'(u)\delta_\star\Delta u
=\underbrace{-F'(u)\mu- \delta^\star F"(u)|\nabla u|^2}_{\leq 0} + \delta_\star\Delta F(u).
$$
The argument can be made rigorous by working on the weak variational formulation of the equation, with suitable approximation of the solution $u$.

\noindent
For proving item (b), we compute
\[ \begin{array}{lll}\ds\frac12\partial_t (v^2\varphi^2)&=&\delta^\star\varphi^2v\nabla\cdot\nabla v-\nu \varphi^2 v
\\
&=&\delta^\star\nabla\cdot(\varphi^2 v\nabla v)-\nu \varphi^2 v-\delta^\star\nabla v\cdot \nabla(\varphi^2 v)
\\
& =& 
\delta^\star\nabla\cdot(\varphi^2 v\nabla v)-\nu \varphi^2 v-\delta^\star|\nabla (\varphi v)|^2
+\delta^\star \ v^2\ |\nabla\varphi|^2.
\end{array}\]
The second and third terms of the right hand side are non positive;
the integral of the  last term is dominated by $\delta_\star\|v\|_{L^\infty(Q)}^2 \|\varphi\|_{H^1(\Omega)}.$
Again a full justification proceeds through an approximation argument.
\QED

For proving Proposition~\ref{osc}, we shall work with several subdomains of $Q_2$, as indicated by Fig.~\ref{Fig:TimeShift}
which might help to follow the arguments.

\begin{lemma}\label{l:3sets}
Let $u$ satisfy $\partial_t u-\delta^\star\Delta u\leq 0$ and $-1\leq u(t,x)\leq +1$ in $Q_2$.
Let us set
\[\begin{array}{l}
\mathscr A=\big\{(t,x)\in Q_1,\ u(t,x)\geq 1/2\big\},
\\
\mathscr B=\big\{(t,x)\in \widetilde Q,\ u(t,x)\leq 0\big\},
\\
\mathscr C=\big\{(t,x)\in Q_1\cup  \widetilde Q,\ 0<u(t,x)<1/2\big\}.\end{array}
\]
There exists $\alpha>0$ such that if
$\mathrm{meas}(\mathscr A)\geq \eta$ and $\mathrm{meas}(\mathscr B)\geq \frac12\mathrm{meas}(\widetilde Q)$, then
$\mathrm{meas}(\mathscr C)\geq \alpha$.
\end{lemma}

\noindent
{\bf Proof.}
We argue by contradiction, assuming that a sequence $\big(u_k\big)_{k\in\mathbb N}$ 
of solutions of 
$\partial_t u_k-\delta^\star\Delta u_k\leq 0$ in $Q_2$
satisfies $-1\leq u_k(t,x)\leq +1$ and 
\begin{equation}\label{hyp_neg}\begin{array}{l}
\mathrm{meas}(\mathscr A_k)\geq \eta,\ \textrm{ with $\mathscr A_k=
\big\{(t,x)\in Q_1,\ u_k(t,x)\geq 1/2\big\}$},
\\[.3cm]
 \mathrm{meas}(\mathscr B_k)\geq \ds\frac12\mathrm{meas}(\widetilde Q),\
 \textrm{ with $\mathscr B_k=
\big\{(t,x)\in \widetilde Q,\ u_k(t,x)\leq 0\big\}$},
\\[.3cm]
\mathrm{meas}(\mathscr C_k)\leq \ds\frac1k,\
  \textrm{ with $\mathscr C_k=
\big\{(t,x)\in Q_1\cup \widetilde Q,\ 0<u_k(t,x)<1/2\big\}$}.
\end{array}\end{equation}
We focus our interest on the positive part $ v_k=[u_k]_+$, with $[z]_+=\max(z,0)$, which is still uniformly bounded: $0\leq v_k(t,x)\leq 1$,
By virtue of Lemma~\ref{distrib}-(a), it satisfies
\begin{equation}\label{eqvk}
\partial_t v_k-\delta^\star \Delta v_k+\mu_k=0,
\end{equation}
with $\mu_k$ a non negative measure.
The strategy can be recapped as follows. We shall establish the compactness of $v_k$ in the 
reduced domain $(-4,0)\times B_{3/2}$. It allows us to assume that $v_k$ converges to a certain function $v$.
Roughly speaking, we are going to show that $v(s,x)$ vanishes on $B_1$ for certain times $-3/2<s<-1$, which
will imply that $v$ vanishes over $Q_1$. It will eventually  lead to a contradiction by considering the behavior of
 the sets $\mathscr A_k$, $\mathscr B_k$, $\mathscr C_k$ as $k\to \infty$.
\\

Let us pick a trial function $\zeta\in C^\infty_c(B_2)$ such that $\zeta(x)=1$ for any $x\in B_{3/2}$ and 
$0\leq \zeta(x)\leq 1$ for any $x\in \mathbb R^N$.
By using Lemma~\ref{distrib}-(b), we get 
for 
$-4<t_1<t_2<0$
\begin{equation}\label{energ_k_2}
\ds\int \zeta^2 |v_k|^2(t_2,x)\ud x +
\delta^\star\ds\int_{t_1}^{t_2}\ds\int  |\nabla (\zeta v_k)|^2(s,x)\ud x\ud s
\leq 
\ds\int \zeta^2 v_k^2(t_1,x)\ud x + C(t_2-t_1),\end{equation}
for a certain constant $C>0$.
In particular, $\big(\zeta v_k\big)_{k\in\mathbb N}$ is bounded in 
$L^\infty(-4,0;L^2(B_2))\cap L^2(-4,0;H^1(B_2))$.
Going back to \eqref{eqvk}, since $\mu_k\geq 0$, $v_k\geq 0$, we observe that
\[\begin{array}{lll}
0\leq 
\ds\int_{t_1}^{t_2}\ds\int_{B_{3/2}}  \mu_k \ud x\ud s
&\leq & \ds\int_{t_1}^{t_2}\ds\int_{B_2} \zeta \mu_k \ud x\ud s
\\[.3cm]
&\leq& 
\ds\int_{B_2} \zeta v_k(t_1,x)\ud x 
- \delta^\star\ds\int_{t_1}^{t_2} \ds\int_{B_2} \nabla v_k\cdot\nabla \zeta\ud x\ud s
\\ [.3cm]
&\leq& 
\|\zeta\|_{L^1}+ 2 \delta^\star \|\nabla v_k\|_{L^2(Q_2)}  \|\nabla\zeta\|_{L^2(B_2)}
\end{array}\]
is bounded uniformly with respect to $k$.
Coming back to \eqref{eqvk}, we deduce that $\big(\partial_t v_k\big)_{k\in\mathbb N}$ is bounded in 
$\mathscr M^1((-4,0)\times B_{3/2})+L^2(-4,0;H^{-1}(B_{3/2}))$.
By virtue of Aubin-Lions-Simon's lemma \cite{JS} (in fact we use the extended version \cite[Theorem~1]{Moussa} which allows us to deal with measure
valued time derivatives), we conclude that 
$\big( v_k\big)_{k\in\mathbb N}$ is compact in $L^2((-4,0)\times B_{3/2})$.
We can thus assume that $v_k$ (possibly relabelling the sequence) converges to some $v$ in $L^2((-4,0)\times B_{3/2})$.
Bienaym\'e-Tchebyschev's inequality yields
\[
\mathrm{meas}\big(\big\{(t,x)\in ((-4,0)\times B_1),\ |v_k(t,x)-v(t,x)|\geq \epsilon\big\}\big)
\leq \ds\frac{\| v_k-v\|^2_{L^2((-4,0)\times B_1)}}{\epsilon^2}\xrightarrow[k\to \infty]{} 0,\] for any $\epsilon>0$.

Let $(t,x)\in (-4,0)\times B_1$ be such that $\epsilon\leq v(t,x)\leq 1/2-\epsilon$.
Then we distinguish the following two cases: either $|v-v_k|(t,x)\geq
 \epsilon$ or 
 $0\leq v_k(t,x)=(v_k-v)(t,x)+v(t,x)\leq |v-v_k|(t,x) +v(t,x)\leq 1/2$.
 It follows that 
 \[\begin{array}{l}
 \mathrm{meas}\big(\big\{(t,x)\in Q_1\cup \widetilde Q,\ \epsilon\leq v(t,x)\leq 1/2-\epsilon  \big\}\big)
 \\[.3cm]
 \qquad
 \leq 
  \mathrm{meas}\big(\big\{(t,x)\in Q_1\cup \widetilde Q,\ |v-v_k|(t,x)\geq\epsilon  \big\}\big)
 \\[.3cm]
 \qquad\qquad + 
   \underbrace{\mathrm{meas}\big(\big\{(t,x)\in Q_1\cup \widetilde Q),\ 0\leq v_k(t,x)\leq1/2  \big\}\big)}_{\mathrm{meas}(\mathscr C_k)}
   \\[.3cm]
 \qquad
 \leq 
  \mathrm{meas}\big(\big\{(t,x)\in Q_1\cup \widetilde Q,\ |v-v_k|(t,x)\geq\epsilon  \big\}\big)
  + 
   \ds\frac 1k,
\end{array}\]
by using \eqref{hyp_neg}.
Letting $k$ go to $\infty$ yields
\[\mathrm{meas}\big(\big\{(t,x)\in Q_1\cup \widetilde Q,\ \epsilon\leq v(t,x)\leq 1/2-\epsilon  \big\}\big)=0.\]
Since this property holds for any $\epsilon$, the monotone convergence property leads to 
\[\mathrm{meas}\big(\big\{(t,x)\in Q_1\cup \widetilde Q,0< v(t,x)< 1/2  \big\}\big)=0.\]
Therefore, we have
\begin{equation}\label{dich}
\textrm{for a.\ e.\ $t\in (-9/4,0)$, either $v(t,x)=0$ or $v(t,x)\geq 1/2$ in $B_1$}.\end{equation}

Similarly, let $(t,x)\in (-4,0)\times B_1$ be such that $ v_k(t,x)=0$.
We distinguish the following two cases: either $|v-v_k|(t,x)\geq
 \epsilon$ or 
 $0\leq v(t,x)=(v-v_k)(t,x)\leq |v-v_k|(t,x) \leq \epsilon$.
Coming back to \eqref{hyp_neg}, we get
\[\begin{array}{l}
\ds\frac 12 \ \mathrm{meas}(\widetilde Q)
\leq \mathrm{meas}(\mathscr B_k)
\\
[.3cm]
\qquad\qquad\qquad
\leq \mathrm{meas}\big(\big\{(t,x)\in \widetilde Q,\ |v-v_k|(t,x)\geq
 \epsilon
\big\}\big)
+
\mathrm{meas}\big(\big\{(t,x)\in \widetilde Q,v(t,x)\leq
 \epsilon
\big\}\big).
\end{array}\]
Letting $k$ go to $\infty$ we obtain 
\[
\ds\frac 12 \ \mathrm{meas}(\widetilde Q)
\leq \mathrm{meas}\big(\big\{(t,x)\in \widetilde Q,v(t,x)\leq
 \epsilon
\big\}\big). \]
By monotone convergence, as $\epsilon\to 0$, we arrive at 
\[
\ds\frac 12 \ \mathrm{meas}(\widetilde Q)
\leq \mathrm{meas}\big(\big\{(t,x)\in \widetilde Q,v(t,x) =0
\big\}\big).
\]
Consequently, we can find a non negligible set of times $s\in (-3/2,-1)$ such that 
$v(s,x)=0$  holds for a.\ e.\ $x\in B_1$.
Letting $k$ go to $\infty$
in \eqref{eqvk}, we obtain $\partial_t v-\delta^\star\Delta v+\nu=0$ on $(-4,0)\times B_{3/2}$,
with $\nu$ a non negative measure.
Let $\zeta\in C^\infty_c(B_{3/2})$ be a non negative trial function such that $\zeta(x)=1$ for any $x\in B_1$.
We apply Lemma~\ref{distrib}-(b), and we obtain for a.~e.~$t\in (s,0)$,
$$
\ds\int_{B_1} v^2(t,x)\ud x\leq \ds\int_{B_{3/2}} v^2(t,x)\zeta^2(x)\ud x\leq  \ds\int_{B_{3/2}} v^2(s,x)\zeta(x)\ud x + C(t-s)=C(t-s),$$
where, owing to \eqref{dich}, we also know that the left hand side is either null or larger than $\frac{\mathrm{meas}(B_1)}{4}$.
We deduce that, actually, $v$ vanishes on $Q_1$. We are going to show that it contradicts \eqref{hyp_neg}.

Indeed, let us consider $(t,x)\in Q_1$ such that  $v_k(t,x)\geq 1/2$. Then, for any $\epsilon>0$, either $|v-v_k|(t,x)
\geq \epsilon$ or $v(t,x)=v_k(t,x)+(v-v_k)(t,x)\geq v_k(t,x) -|v-v_k|(t,x) \geq 1/2-\epsilon$.
With the first property in \eqref{hyp_neg}, it follows that 
\[\begin{array}{lll}
\eta\leq \mathrm{meas}(\mathscr A_k)&\leq& 
\mathrm{meas}\big(\big\{(t,x)\in Q_1,\ |v-v_k|(t,x)
\geq \epsilon\big\}\big)
\\[.3cm]
&&\qquad
+ \mathrm{meas}\big(\big\{(t,x)\in Q_1,v(t,x)
\geq 1/2-\epsilon\big\}\big).
\end{array}\]
Letting $k$ go to $\infty$ yields
\[\eta\leq \mathrm{meas}\big(\big\{(t,x)\in Q_1,v(t,x)
\geq 1/2-\epsilon\big\}\big).\]
Since this inequality holds for any $\epsilon>0$, we conclude, by monotone convergence, that 
\[\eta\leq 
 \mathrm{meas}\big(\big\{(t,x)\in Q_1,v(t,x)
\geq 1/2\big\}\big)
\]
holds, a contradiction.
 \QED

\noindent
{\bf Proof of Proposition~\ref{osc}.}
We consider $(t,x)\mapsto v(t,x)$ such that $-1\leq v(t,x)\leq +1$, $\mathrm{meas}\big(\big\{ 
(t,x)\in \widetilde Q,\ v(t,x)\leq 0\big\}\big)\geq \mu\ \mathrm{meas}(\widetilde Q)$, and $v$ satisfies $\partial_t v-\delta^\star \Delta v\leq 0$ in $Q_2$.
The proof splits into two steps.

\noindent
{\it \underline{Step 1.}}

\noindent
For $k\in\mathbb N$, set 
$$v_k(t,x)=2^k(v(t,x)-(1-1/2^k)).$$
We shall show that  the integral
\[
\ds\iint_{Q_1} [v_k]^2_+\ud x\ud t\]
can be made as small as we wish, by choosing $k$ large enough. 
Observe that
$$v_k=2^k(v-1)+1=2v_{k-1}-1$$
which implies that $v_k\leq 1$ and
$$\big\{(t,x)\in \widetilde Q,\ v(t,x)\leq 0\big\}\subset \big\{(t,x)\in \widetilde Q,\ v(t,x)\leq 1-1/2^k\big\}=
\big\{(t,x)\in \widetilde Q,\ v_k(t,x)\leq 0\big\}.$$
Thus, by assumption on $v$,we have
$$
\mathrm{meas}\big(\big\{ (t,x)\in \widetilde Q,\ v_k(t,x)\leq 0\big\}\big)\geq 
\mathrm{meas}\big(\big\{ (t,x)\in \widetilde Q,\ v(t,x)\leq 0\big\}\big)\geq \mu\ \mathrm{meas}(\widetilde Q).$$
Let us suppose that, for any $k\in\mathbb N$ 
\[
\ds\iint_{Q_1} [v_k]_+^2\ud x\ud t\geq \delta
 \]
 holds for a certain $\delta >0$.
Since this integral is dominated by 
$$
\mathrm{meas}\big(\big\{ (t,x)\in  Q_1,\ v_k(t,x)\geq 0\big\}\big)=
\mathrm{meas}\big(\big\{ (t,x)\in  Q_1,\ v_{k-1}(t,x)\geq 1/2\big\}\big)
$$
we infer
\[\mathrm{meas}\big(\big\{ (t,x)\in  Q_1,\ v_{k-1}(t,x)\geq 1/2\big\}\big)\geq \delta\]
independently of $k$.
Applying Lemma~\ref{l:3sets} yields
\[\mathrm{meas}\big(\big\{ (t,x)\in  Q_1\cup \widetilde Q,\ 0<v_{k-1}(t,x)< 1/2\big\}\big)\geq \alpha,\]
still independently of $k$.
It follows that 
\[\begin{array}{l}
\mathrm{meas}\big(\big\{ (t,x)\in  Q_1\cup \widetilde Q,\ v_{k}(t,x)\leq 0\big\}\big)
\\[.3cm]
\qquad
=\mathrm{meas}\big(\big\{ (t,x)\in  Q_1\cup \widetilde Q,\ 2v_{k-1}(t,x)-1\leq 0\big\}\big)
\\[.3cm]
\qquad
=\mathrm{meas}\big(\big\{ (t,x)\in  Q_1\cup \widetilde Q,\ v_{k-1}(t,x)\leq 0\big\}\big)
\\
[.3cm]
\qquad\qquad+\mathrm{meas}\big(\big\{ (t,x)\in  Q_1\cup \widetilde Q,\ 0<v_{k-1}(t,x)\leq  1/2\big\}\big)
\\[.3cm]
\qquad 
\geq 
\mathrm{meas}\big(\big\{ (t,x)\in  Q_1\cup \widetilde Q,\ v_{k-1}(t,x)\leq 0\big\}\big)+\alpha.
\end{array}\]
Since $\mathrm{meas}\big(\big\{ (t,x)\in  Q_1\cup \widetilde Q,\ v_{0}(t,x)\leq 0\big\}\big)
\geq \mathrm{meas}\big(\big\{ (t,x)\in   \widetilde Q,\ v_{0}(t,x)\leq 0\big\}\big)\geq \mu\ \mathrm{meas}(\widetilde Q)$,
this recursion formula leads to 
\[
\mathrm{meas}\big(\big\{ (t,x)\in  Q_1\cup \widetilde Q,\ v_{k}(t,x)\leq 0\big\}\big)\geq 
 \mu\ \mathrm{meas}(\widetilde Q)+k\alpha
.
\]
However, this cannot occur for any $k$ since the left hand side is bounded by $\mathrm{meas}(Q_2)$.
We conclude that, given $\delta >0$, there exists $k_\star\in\mathbb N$ such that 
 \[\ds\iint_{Q_1} [v_{k_\star}]_+^2\ud x\ud t\leq \delta.\]
 \\
 
 \noindent
 {\it\underline {Step 2.}}
 
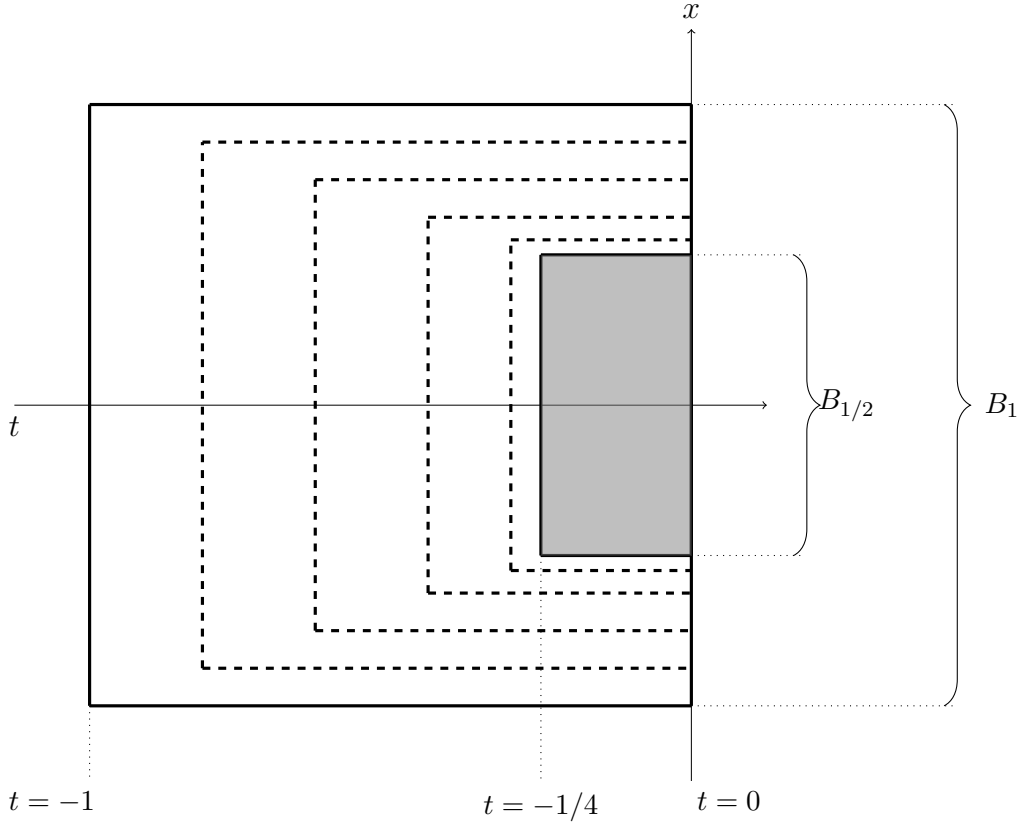
\begin{figure}[!ht]
\begin{center}
\begin{tikzpicture}
\draw [very thin,->]  (0,-5) -- (0,5);
\draw (0,5) node[above] {$x$} ;
\draw [very thin,->]  (-9,0) -- (1,0);
  \draw (-9,0) node[below] {$t$};

   \draw [decorate,decoration={brace,amplitude=10pt},xshift=-4pt,yshift=0pt]
(1.5,2) -- (1.5,-2)node [black,midway,xshift=20pt] {\small $B_{1/2}$};

 \draw [decorate,decoration={brace,amplitude=10pt},xshift=-4pt,yshift=0pt]
(3.5,4) -- (3.5,-4)node [black,midway,xshift=20pt] {\small $\ B_{1}$};
  
  \draw [very thick]  (-8,-4) -- (-8,4) ;
  \draw [very thick]  (-8,4) -- (0,4) ;
   \draw [very thick]  (0,4) -- (0,-4) ;
   \draw [very thick]  (-8,-4) -- (0,-4) ;
  
  \draw [dotted] (-8,-4) -- (-8,-5);
    
       \draw [dashed, very thick] (-6.5,-3.5) -- (0,-3.5);
    \draw [dashed, very thick] (-6.5,3.5) -- (0,3.5);
     \draw [dashed, very thick] (-6.5,-3.5) -- (-6.5,3.5);

    \draw [dashed, very thick] (-5,-3) -- (0,-3);
    \draw [dashed, very thick] (-5,3) -- (0,3);
     \draw [dashed, very thick] (-5,-3) -- (-5,3);
     
  \draw [dashed, very thick] (-3.5,-2.5) -- (0,-2.5);
    \draw [dashed, very thick] (-3.5,2.5) -- (0,2.5);
     \draw [dashed, very thick] (-3.5,-2.5) -- (-3.5,2.5);
     
       \draw [dashed, very thick] (-2.4,-2.2) -- (0,-2.2);
    \draw [dashed, very thick] (-2.4,2.2) -- (0,2.2);
     \draw [dashed, very thick] (-2.4,-2.2) -- (-2.4,2.2);
    
    \draw [dotted] (0,-4) -- (3.5,-4);
   \draw [dotted] (0,4) -- (3.5,4);
    \draw [dotted] (0,-2) -- (1.5,-2);
     \draw [dotted] (0,2) -- (1.5,2);
      \draw [dotted] (-2,-2) -- (-2,-5);
  
   \draw (-8.5,-5) node[below] {{\small  {$t=-1$}}};
  \draw (-2,-5) node[below] {{\small  {$t=-1/4$}}};
  \draw (0.5,-5) node[below] {{\small  {$t=0$}}};
  
  \draw [very thick]  (-2,-2) -- (-2,2) ;
  \draw [very thick]  (-2,2) -- (0,2) ;
   \draw [very thick]  (0,2) -- (0,-2) ;
   \draw [very thick]  (-2,-2) -- (0,-2) ;


  
\coordinate[] (A) at (-2,-2);
\coordinate[] (B) at (0,-2);
\coordinate[] (C) at (0,2);
\coordinate[] (D) at (-2,2);
%
%
%
%
\fill[gray, opacity=0.5, draw=black] 
(A) -- (B) -- (C) -- (D) -- cycle;


%
%
      
 \end{tikzpicture}
\caption{
The domains $Q_1$,  $Q_{r_\ell}$ and $Q_{1/2}$ (the grey box)}
\label{Fig:Qell}\end{center}
\end{figure}

  \noindent
  The second step relies on  De Giorgi's analysis.
  Let us set $w(t,x)=v_{k_\star}(t,x)$.
  We shall show that, provided $\delta$ is small enough (which means $k_\star$ large enough),
  $w(t,x)\leq 1/2$ on $Q_{1/2}$.
  To this end, let us set, for $\ell\in \mathbb N$,
  \[\begin{array}{l}
   m_\ell=\ds\frac12\Big(1-\ds\frac{1}{2^\ell}\Big),
  \\[.3cm]
  w_\ell(t,x)=[w(t,x)-m_\ell]_+,
  \\[.3cm]
  r_\ell= \ds\frac12\Big(1+\ds\frac{1}{2^\ell}\Big), \qquad 
  t_\ell= -r_\ell^2=- \ds\frac14\Big(1+\ds\frac{1}{2^\ell}\Big)^2.
  \end{array}\] 
  We are going to work in the domains $Q_{1/2}\subset Q_{r_\ell}\subset Q_1$ that shrink to $Q_{1/2} $ as $\ell\to \infty$,
   see Fig.~\ref{Fig:Qell}.

We consider a sequence of functions $\zeta_\ell \in C_c^\infty(B_{r_{\ell-1}})$
such that $0\leq \zeta_\ell(x)\leq 1$ on $B_{r_{\ell-1}}$ and $\zeta_\ell(x)= 1$ on $B_{r_{\ell}}$.
We shall use the basic estimate
\[
|\nabla\zeta_\ell(x)|\leq C2^\ell,\qquad \ds\frac{1}{t_\ell-t_{\ell-1}}\leq C2^{2\ell}.\]
We already know that $0\leq w_\ell(t,x)\leq 1$, by definition. We can apply the energy estimate in Lemma ~\ref{distrib}, which reads
\begin{equation}\label{energ00}
\begin{array}{l}
 \ds\frac12\ds\int_{B_1} w_\ell^2(t,x)\zeta_\ell^2(x)\ud x 
 +\delta^\star\ds\int_s^t\ds\int_{B_1}| \nabla(\zeta_\ell w_\ell)|^2(\tau, x)\ud x\ud \tau
\\
[.3cm]
\qquad\qquad\qquad\qquad \leq  \ds\frac12\ds\int_{B_1} w_\ell^2(s,x)\zeta_\ell^2(x)\ud x 
+\delta^\star
\ds\int_s^t\ds\int_{B_1} w_\ell^2| \nabla\zeta_\ell |^2 (\tau, x)\ud x\ud \tau.
 \end{array}\end{equation} 
 for $-1<s< t_\ell < t<0$
(note that here we keep explicit the integral in the right hand side
 that is roughly estimated by a constant in Lemma~\ref{distrib}).
Averaging over 
$s\in (t_{\ell-1},t_\ell)$ (and using the fact that the integral of a positive quantity 
over $(s,t)$ is thus bounded below --- resp. above --- by the  
integral over $(t_\ell,t)$ --- resp. $(t_{\ell-1},t)$)
yields
$$
\begin{array}{l}
 \ds\frac12\ds\int_{B_1} w_\ell^2(t,x)\zeta_\ell^2(x)\ud x 
 +\delta^\star\ds\int_{t_{\ell}}^t\ds\int_{B_1}| \nabla(\zeta_\ell w_\ell)|^2(\tau, x)\ud x\ud \tau
\\
[.3cm]
\qquad\qquad\qquad\qquad \leq (1/2+\delta^\star)C 2^{2\ell}
\ds\int_{t_{\ell-1}}^0 \ds\int_{\mathrm{supp}(\zeta_\ell)} |w_\ell|^2(\tau ,x) \ud x\ud \tau.
 \end{array}$$
 Let us set 
 \[\mathscr U_\ell=\ds\int_{t_\ell}^0\ds\int_{B_{r_\ell}} |w_\ell|^2(t,x)\ud x\ud t,\]
 and
  \[\mathscr E_\ell=\ds\sup_{t_\ell\leq t\leq 0}\ds\int_{B_1} w_\ell^2(t,x)\zeta_\ell^2(x)\ud x 
 +\ds\int_{t_{\ell}}^0\ds\int_{B_1}| \nabla(\zeta_\ell w_\ell)|^2(\tau, x)\ud x\ud \tau
\]
 We wish to establish a non linear recursion for $\mathscr U_\ell$, which will allow us to justify that it tends to 0 as $\ell\to \infty$.
 On the one hand, since 
 \[w_\ell\leq w_{\ell-1} \qquad\textrm{and} \qquad \mathrm{supp}(\zeta_\ell)\subset B_{\ell-1},\]
 we note that \eqref{energ00} yields 
 \[
 \mathscr E_\ell
  \leq 
   (2+1/\delta^\star) (1/2+\delta^\star)C 2^{2\ell}
 \mathscr U_{\ell-1}.
 \]
 On the other hand, we observe that
$$\begin{array}{lll} 
 \mathscr U_\ell&\leq& 
 \ds\int_{t_\ell}^0\ds\int_{B_{r_\ell}} |\zeta_\ell w_\ell|^2(t,x)\ud x\ud t
 \\
 [.3cm]&\leq&
 \left(  \ds\int_{t_\ell}^0\ds\int_{B_{r_\ell}} |\zeta_\ell w_\ell|^{2(N+2)/N}(t,x)\ud x\ud t \right)^{N/(N+2)}\\
 [.3cm]&&
 \qquad\qquad\qquad \times
 \left(\mathrm{meas}\big(\big\{(t,x)\in (t_\ell,0)\times B_{r_\ell},\ \zeta_\ell w_\ell>0\big\}\big)  \right)^{2/(N+2)}
  ,\end{array}$$
  by using H\"older's inequality.
  Remark that 
  \[w-m_{\ell-1}=w-m_\ell+\ds\frac{1}{2^{\ell+1}}
  \]
  which leads to to
  \[\begin{array}{l}
  \mathrm{meas}\big(\big\{(t,x)\in (t_\ell,0)\times B_{r_\ell},\ \zeta_\ell w_\ell>0\big\}\big) 
  \\[.3cm]
  \qquad
 \leq   \mathrm{meas}\big(\big\{(t,x)\in (t_{\ell-1},0)\times B_{r_{\ell-1}},\  w_{\ell-1}>2^{-\ell-1}\big\}\big) 
\\[.3cm]
  \qquad
 \leq 2^{2\ell+2} \mathscr U_{\ell-1},
  \end{array}\]
    by virtue of the Bienaym\'e-Tchebyschev inequality.
    Next, we use the  
    Gagliardo-Nirenberg-Sobolev inequality, see \cite[Theorem p.~125]{Nir}
    \[ \left( \ds\int_{B_{r_\ell}}   |\zeta_\ell w_\ell|^{2N/(N-2)}(t,x)\ud x\right)^{(N-2)/N}
  \leq C_S
   \ds\int_{B_{r_\ell}}   |\nabla(\zeta_\ell w_\ell)|^{2}(t,x)\ud x.\]
    Mind that we have integrated with respect to the space variable only.
   We can write
   \[
\ds\frac{N+2}{N}=\theta\ds\frac{2N}{N-2}+2(1-\theta),
\qquad
\theta=\ds\frac{N-2}{N}\in(0,1),\]
so that
\[\begin{array}{l}
\ds\int_{t_\ell}^0 \ds\int_{B_{r_\ell}}   |\zeta_\ell w_\ell|^{2(N+2)/N}(t,x)\ud x\ud t
\\[.3cm]
\qquad\leq
\ds\int_{t_\ell}^0 
\left(
 \ds\int_{B_{r_\ell}}   |\zeta_\ell w_\ell|^{2N/(N-2)}(t,x)\ud x
\right)^{\theta}
\underbrace{ \left( \ds\int_{B_{r_\ell}}   |\zeta_\ell w_\ell|^{2}(t,x)
\ud x
\right)^{1-\theta}}_{\leq  \mathscr E_\ell^{1-\theta}}
\ud t
\\[.3cm]
\qquad \leq C_S^\theta \mathscr E_\ell^{2-\theta}.
\end{array}
\]
 Therefore, gathering all these informations together,  we obtain
\[
\mathscr U_\ell\leq \Lambda ^\ell
 \mathscr U_{\ell-1}^{1+2/(N+2)},
 \]
 for a certain constant $\Lambda >1$.
Owing to Lemma~\ref{2bete}, we deduce that $\lim_{\ell\to \infty} \mathscr U_\ell=0$ provided $\mathscr U_0$ is small enough.
The smallness condition on $\mathscr U_0$ is precisely ensured by the definition  $w=v_{k_\star}$ coming from Step~1.
Since 
$$\ds\frac{1}{|t_{\ell}|}\ds\int_{t_\ell}^0\ds\int_{B_{r_\ell}} |w_\ell|^2(t,x)\ud x\ud t\leq \mathscr U_\ell$$
we conclude, by applying Fatou's lemma, that
$$ 2\ds\iint_{Q_{1/2}} [w-1/2]_+^2(t,x)\ud x\ud t\leq \liminf_{\ell\to\infty}\ds\frac{1}{t_{\ell}}\ds\int_{t_\ell}^0\ds\int_{B_{r_\ell}} |w_\ell|^2(t,x)\ud x\ud t=0$$
so that, finally, $w(t,x)\leq 1/2$ holds a.~e.~on $Q_{1/2}$.

Coming back to the change of unknown $w(t,x)=v_{k_\star}(t,x)=2^{k_\star}(v(t,x)-(1-1/2^{k_\star}))\leq 1/2$
becomes $$
v(t,x)\leq 1+\ds\frac{1}{2^{k_\star+1}}-\ds\frac{1}{2^{k_\star}}=1-\ds\frac{1}{2^{k_\star+1}}<1.
$$
 \QED

\section{$L^{(N+1)/N}$ estimate on the total mass}
\label{S:4}

This Section is devoted to the proof of the following statement.

\begin{proposition}\label{Fabes}
There exists a constant $K>0$ such that, $M\geq 0$ being a solution of \eqref{eqMass} in $Q_2$. We have
\[
\|M\|_{L^{(N+1)/N}(Q_1)}\leq K \ \ds\sup_{-4\leq t\leq 0} \ds\int_{B_2} M(t,x)\ud x
.\]
\end{proposition}

\noindent
{\bf Proof.}
Let $f$ be in $C^\infty_c(Q_1)$ such that 
$$\|f\|_{L^{N+1}(Q_1)}\leq 1.$$
We consider the solution of the final problem 
\begin{equation}\label{equ}\begin{array}{ll}
\partial_t u +d\Delta u=f \qquad & \textrm{in $(0,T)\times \mathbb R^N$},
\\
u(T,x)=0,\qquad&  u\big|_{\partial B_2}=0.
\end{array}\end{equation}
We start by reminding the reader the Alexandrof-Bakelman-Pucci-Krylov-Tso  (ABPKT) inequality 
\cite{Ale,Bak,Puc,Kry,Tso}: there exists a constant $\mathscr C>0$ such that
\begin{equation}\label{ABP}
\ds\sup_{(t,x)\in Q_2}|u(t,x)| \leq \mathscr C\  \|f\|_{L^{N+1}(Q_2)}.
\end{equation}
In order to obtain an estimate on the $L^{(N+1)/N}(Q_1)$ norm of $M$, solution of \eqref{eqMass}, 
we proceed by duality, bearing in mind the definition
\[
\|M\|_{L^{(N+1)/N}(Q_1)}=
\sup \left\{
\Big|\ds\iint_{Q_1} Mf\ud x\ud t\Big|,\ 
f\in C^\infty_c(Q_1),\ \|f\|_{L^{N+1}(Q_1)}\leq 1\right\}.\]
\\

Let $\zeta$ be a cut-off function: 
$\zeta\in C^\infty_c(B_{3/2})$, $\zeta(x)=1$ for any $x\in B_1$, and $0\leq \zeta(x)\leq 1$ for any $x\in\mathbb R^N$.
Remark that
\[
\ds\iint_{Q_2} \zeta  M f\ud x\ud t=\ds\iint_{Q_1} Mf \ud x\ud t,\]
since $\mathrm{supp}(f)\subset Q_1$.
We compute this integral by using \eqref{equ}
\[
\begin{array}{lll}
\ds\iint_{Q_2} \zeta  M f\ud x\ud t
&=&
\ds\iint_{Q_2} \zeta  M(\partial_t u+d\Delta u) \ud x\ud t
\\[.3cm]
&=&\ds\int_{-2}^0 \ds\frac{\ud}{\ud t}\left( \ds\int_{B_2} \zeta  Mu \ud x\right)\ud t
-\ds\iint_{Q_2} \zeta  u\Delta (dM) \ud x\ud t
\\[.3cm]&&+\ds\iint_{Q_2} \zeta  Md\Delta u \ud x\ud t
\\[.3cm]
&=& \ds\int_{B_2} \zeta  Mu(0,x) \ud x
-2\ds\iint_{Q_2}dM  \nabla\zeta \cdot \nabla u\ \ud x\ud t
\\[.3cm]&&-\ds\iint_{Q_2} udM \Delta \zeta  \ud x\ud t.
\end{array}\]
We have used several integration by parts where the boundary terms vanish owing to the fact that 
$\mathrm{supp}(\zeta)\subset B_{3/2}\subset B_2$.
The integrand of the penultimate in the right hand side 
can be rewritten as $\sqrt {dM}\nabla u\cdot  \sqrt {dM}\nabla\zeta$, 
and then we use the Cauchy-Schwarz inequality and the Young inequality 
$ab=\sqrt{\kappa} a\frac{b}{\sqrt {\kappa}}\leq \frac12(\kappa a^2+ \frac{b^2}{\kappa})$.
We thus arrive at the following estimate
\begin{equation}\label{est0}
\begin{array}{lll}
\left|\ds\iint_{Q_2} \zeta f M\ud x\ud t\right|
&\leq& \left| \ds\int_{B_2} \zeta  Mu(0,x) \ud x\right|
\\[.3cm]&&+ 
\kappa \ds\iint_{Q_2}dM   | \nabla u|^2\ \ud x\ud t
+\ds\frac 1\kappa  \ds\iint_{Q_2}dM   | \nabla \zeta|^2\ \ud x\ud t
\\[.3cm]&&+\left|\ds\iint_{Q_2} udM \Delta \zeta   \ud x\ud t\right|,
\end{array}\end{equation}
where $\kappa\in(0,1)$ is a parameter that will be
determined later on.
Inspired from \cite[proof of Theorem~2.1]{FaSt}, in order to estimate the second integral in the right hand side, we use the elementary relation 
\[
|\nabla u|^2=\ds\frac12 \Delta (u^2)-u\Delta u  .\]
Going back to  \eqref{equ}, we are thus led to 
\[d|\nabla u|^2=\ds\frac d 2 \Delta (u^2)+ \ds\frac12\partial_t (u^2)  - uf.\]
The advantage of this formulation relies on the fact 
that, denoting $\nu$ the outward unit normal on $\partial B_2$, 
\[u\big|_{\partial B_2}=u^2\big|_{\partial B_2}=0,\qquad
\nabla u^2\cdot \nu\big|_{\partial B_2} =2u\nabla u\cdot \nu\big|_{\partial B_2}=0,\]
which allows us to perform further integration by parts. 
We get
\[\begin{array}{l} \ds\iint_{Q_2}dM   | \nabla u|^2\ \ud x\ud t
\\[.3cm]
\qquad
= \ds\frac12\ds\iint_{Q_2}dM  \Delta ( u^2)\ \ud x\ud t + \ds\frac12 \ds\iint_{Q_2} M  \partial_t  (u^2)\ \ud x\ud t
- \ds\iint_{Q_2}M  u f \ud x\ud t
\\[.3cm]
\qquad
= -\ds\frac12\ds\iint_{Q_2}\nabla(dM)\cdot  \nabla ( u^2)\ \ud x\ud t +  \ds\frac12 \ds\int_{B_2} M   u^2(0,x)\ \ud x
\\[.3cm]
\qquad\qquad- \ds\frac12 \ds\iint_{Q_2}\Delta(dM)   u^2 \ud x\ud t
-  \ds\iint_{Q_2}M  u f \ud x\ud t
\\[.3cm]
\qquad
= \ds\frac12\ds\int_{B_2}M   u^2(0,x)\ \ud x
-  \ds\iint_{Q_2}M  u f \ud x\ud t.
\end{array}\]
For the last term, since $\mathrm{supp}(f)\subset Q_1$, the integral actually reduces over $Q_1$ only.
The H\"older inequality then yields
\[\begin{array}{lll}
 \left| \ds\iint_{Q_2}M  u f \ud x\ud t\right|
 &=& \left| \ds\iint_{Q_1}M  u f \ud x\ud t\right|
 \leq \|u\|_{L^\infty(Q_1)}  \|M\|_{L^{(N+1)/N}(Q_1)}  \|f\|_{L^{N+1}(Q_1)} 
  \\[.3cm]&\leq&\mathscr C\  \|f\|^2_{L^{N+1}(Q_1)}   \|M\|_{L^{(N+1)/N}(Q_1)}, 
  \end{array}\]
by using \eqref{ABP}.
Besides, still by using \eqref{ABP} and  $\mathrm{supp}(f)\subset Q_1$, we get 
\[\begin{array}{lll}
\ds\frac12\ds\iint_{Q_2}M   u^2(0,x) \ud x &\leq& \ds\frac12 \|u\|_{L^\infty(Q_2)}^2\|M\|_{L^\infty(-4,0;L^1(Q_2))}
\\[.3cm]
& \leq&\mathscr C^2\  \|f\|^2_{L^{N+1}(Q_1)}   \|M\|_{L^\infty(-4,0;L^1(Q_2))}.\end{array}\]
The last two terms in the right hand side of \eqref{est0} are estimated as follows: we get
 \[
  \ds\iint_{Q_2}dM   | \nabla \zeta|^2\ \ud x\ud t
  \leq 4\delta^\star \|\zeta\|^2_{W^{1,\infty}(B_2)}\|M\|_{L^\infty((-4,0);L^1(B_2))},\]
 and
\[
\begin{array}{lll}
\left|\ds\iint_{Q_2} udM \Delta \zeta   \ud x\ud t\right|
&\leq& 4\delta^\star \|\zeta\|_{W^{2,\infty}(B_2)}\   \|u\|_{L^\infty(Q_2)} \|M\|_{L^\infty((-4,0);L^1(B_2))}
\\[.3cm] &\leq& 4\delta^\star  \|\zeta\|_{W^{2,\infty}(B_2)}\ \mathscr C \|f\|_{L^{N+1}(Q_1)} \|M\|_{L^\infty((-4,0);L^1(B_2))}.
\end{array}\]
The first integral in the right hand side of \eqref{est0}
is dominated by
\[
\|u\|_{L^\infty(Q_2)}\|M\|_{L^\infty(-4,0;L^1(Q_2))}
\leq \mathscr C \|f\|_{L^{N+1}(Q_2)} \|M\|_{L^\infty(-4,0;L^1(Q_2))}.
\]
Finally, we have found a constant $C>0$ such that  for any $f\in C^\infty_c(Q_1)$, with $\|f\|_{L^{N+1}(Q_1)}\leq 1$, we have
\[
\left|\ds\iint_{Q_1} fM\ud x\ud t
\right|
\leq C\Big(\Big(1+\kappa + \ds\frac1\kappa\Big) \|M\|_{L^\infty(-4,0;L^1(B_2))} + \kappa  \|M\|_{L^{(N+1)/N}(Q_1)} \Big).\]
Taking the supremum over such $f$'s makes the dual norm $L^{(N+1)/N}(Q_1)$ appear.
We choose $\kappa$ small enough, so that $1-\kappa C>1$, and we conclude that
\[
 \|M\|_{L^{(N+1)/N}(Q_1)}\leq \ds\frac{C(1+\kappa +1/\kappa)}{1-\kappa C} \|M\|_{L^\infty(-4,0;L^1(B_2))}
\]
holds.
\QED

\section{End of proof of Theorem \ref{Theo_principal}:  proof of Lemma~\ref{lemm_de giorgi2}}
\label{S:fin}

\noindent
 Let $0<\epsilon_0<\sqrt{T_{\mathrm{max}}/2}$.  
For each component $a^{(\epsilon)}_i$, Proposition~\ref{Fabes} gives
\begin{equation}\label{lemmadeg1}
\|a^{(\epsilon)}_i\|_{L^{(N+1)/N}(Q_1)}\leq \|M^{(\epsilon)}\|_{L^{(N+1)/N}(Q_1)}\leq K \ \|M^{(\epsilon)}\|_{L^\infty(-4,0;L^1(B_2))}. 
\end{equation}
Next, Lemma~\ref{l:L1shr}, yields 
\begin{equation}\label{lemmadeg2}
\|M^{(\epsilon)}\|_{L^\infty(-4,0;L^1(B_2))}\leq c\|\Phi\|_{L^\infty} \epsilon^{\alpha-2+2/(q-1)}
.
\end{equation}
Combining~(\ref{lemmadeg1}) and ~(\ref{lemmadeg2}) with Proposition~\ref{P:CaVa} leads to
\begin{equation}\label{lemmadeg3}
\sum_{i=1}^{p} \|a^{(\epsilon)}_i\|_{L^{(N+1)/N}(Q_1)}\leq \mathscr K\ \|a^0\|_{L^\infty(\mathbb R^N)}^{1-2/N} \|a^0\|_{L^1(\mathbb R^N)}^{2/N}\ \epsilon^{\alpha-2+2/(q-1)}
\end{equation}
for a constant $\mathscr K$ which depends on $p$ and  $N$.
This information is useful as far as the degree of non linearities is such that the exponent remains positive, which means
  $q\leq 2+\frac{\alpha}{2-\alpha}$.
It ends the proof of Lemma~\ref{lemm_de giorgi2}.

As explained in Section~\ref{S:main}, having at hand this property of the rescaled solution we go back to the original unknown, and we deduce the $L^\infty$ bound of the solution, see Corollary~\ref{unibound}.  Theorem~\ref{Theo_principal}, and therefore Theorem~\ref{theo_main} too, is fully justified.
\QED

\begin{rmk}
The estimates discussed above differ from \cite[see sp. Corollary~14 \& Lemma~15]{CaVa};
and in particular the smallness condition  on $\epsilon_0$ does not involve the initial entropy \eqref{hyp_ci}.
\end{rmk}

\appendix

\section{Proof of~Lemma~\ref{lemm_de giorgi}}
\label{a:L3}

The proof is based on the De Giorgi techniques~\cite{DeG} and it is reminiscent of the method introduced by Alikakos \cite{Alikakos}. 
We exploit the dissipative properties of the system by considering the following 
 non negative,  non decreasing, convex,  and $C^1$ function
\[H(z)=\left\{\begin{array}{ll}
(1+z)\ln(1+z)-z\qquad &\text{ if } z\geq 0,\\
0 \qquad&\text{ if } z\leq 0.\end{array}\right.\]
 Let us introduce  the following  sequences, for  $j\in\mathbb N$,
 $$k_j=1-2^{-j},\qquad t_j=1/4+2^{-j-2}.$$ Henceforth, we set $\mathcal{B}_j=B_{t_j}$
  and $\mathcal{Q}_j=(-t_j,0)\times \mathcal{B}_j$. Note that  $$\begin{array}{l} B\Big(0,\ds\frac{1}{4}\Big)\subset\mathcal{B}_j\subset  \mathcal{B}_{j-1}\subset B\Big(0,\ds\frac{1}{2}\Big),
\\[.3cm]
\Big(-\ds\frac14,0\Big)\times B\Big(0,\ds\frac{1}{4}\Big)\subset \mathcal{Q}_j\subset
  \mathcal{Q}_{j-1}\subset \Big(-\ds\frac12,0\Big)\times B\Big(0,\ds\frac{1}{2}\Big).\end{array}$$
We also introduce a family of  cut-off functions that satisfies the following properties
\[\begin{array}{ll}
\zeta_j:\mathbb R^N\rightarrow [0,\infty),\qquad& \zeta_j\in C^\infty_c(\mathbb R^N),\\
0\leq \zeta_j(x)\leq 1,\qquad &\\
\zeta_j(x)=1 \textrm{ for $x\in \mathcal B_j$},\qquad&
 \zeta_j(x)=0 \textrm{ for $x\in\R^N\setminus \mathcal B_{j-1}$},
 \end{array}\]
 and
 \[
 \ds\sup_{l,m\in \{1,..., N\},\ x\in\mathbb R^N} | \partial_{l,m}^2\zeta_j(x)|\leq C\  2^{2j}\,\,\,\hbox{for a certain constant $C>0$}.
 \]

\begin{lemma}\label{gainofr}
%
There exists a constant $\hat C>0$, which depends only on $\delta_\star$, $\delta^\star$, and on h1)-h4), such that
for any solution $a=(a_1,...,a_p)$ of \eqref{eqA2} and any $\eta\in [0,1]$, we have 
\[\begin{array}{l}
\ds\sup_{-t_j\leq t\leq 0}\ds\sum_{i=1}^{p}\int_{\mathcal B_{j}}
 H(a_i-\eta)(t,x)\ud x
 +
 4\delta_\star
 \ds\sum_{i=1}^{p}\iint_{\mathcal Q_{j}}\big|\nabla_x \sqrt{1+[a_i-\eta]_+}\big|^2(\tau,x)\ud x\ud \tau
\\[.3cm]
\leq 
\hat C\left(
2^{2j}
\ds\sum_{i=1}^{p}\ds\int_{-t_{j-1}}^{0} \int_{\mathcal B_{j-1}}
 H(a_i-\eta)(s,x)\ud x
\ud s\right.
\\
\qquad\qquad\left.
+\ds\sum_{i=1}^{p}
\ds\int_{-t_{j-1}}^0
\int_{\mathcal B_{j-1}}(1+[a_i-\eta]_+)^{q-1}\ln (1+[a_i-\eta]_+)(\tau,x)\ud x\ud \tau\right).
\end{array}\]
\end{lemma}

\noindent
{\bf Proof.}
Multiply \eqref{eqA2} by $\zeta_jH'(a_i-\eta)$,  integrate over ${\mathcal B_{j-1}}$ and sum. We get
\begin{equation}\label{int1}
\begin{array}{l}
\ds\frac{\ud}{\ud t}
\ds\sum_{i=1}^{p}\int_{\mathcal B_{j-1}}
 \zeta_jH(a_i-\eta)\ud x
 \\
 \qquad=\ds \sum_{i=1}^{p} \int_{\mathcal B_{j-1}}
 d_i\Delta a_i\
 H'(a_i-\eta)\zeta_j\ud x+ \ds\sum_{i=1}^{p}\int_{\mathcal B_{j-1}}
 Q_i(a)H'(a_i-\eta)\zeta_j\ud x.
 \end{array}\end{equation}
The first term in the right hand side  of~(\ref{int1}) can be written as 
$$-\ds \sum_{i=1}^{p} \int_{\mathcal B_{j-1}}
 d_i|\nabla a_i|^2\
 H''(a_i-\eta)\zeta_j\ud x+
 \sum_{i=1}^{p} \int_{\mathcal B_{j-1}}
 d_i H(a_i-\eta)\Delta\zeta_j\ud x
$$
where, on the one hand,
 \[\begin{array}{l}
\ds\sum_{i=1}^{p}\int_{\mathcal B_{j-1}}
 d_i|\nabla_x a_i|^2\
 H''(a_i-\eta) \zeta_j\ud x \geq
 4\delta_\star
 \ds\sum_{i=1}^{p}\int_{\mathcal B_{j}}\big|\nabla_x \sqrt{1+[a_i-\eta]_+}\big|^2\ud x,
 \end{array}\]
 and, on the other hand, 
 \[
  \sum_{i=1}^{p} \int_{\mathcal B_{j-1}}
 d_i H(a_i-\eta)\Delta\zeta_j\ud x
 \leq C\delta^\star 2^{2j} \int_{\mathcal B_{j-1}}
H(a_i-\eta)\ud x.
 \]
 For the second term  in the right hand side  of~(\ref{int1}), we get 
 \[
 \begin{array}{l}
\ds\sum_{i=1}^{p}\int_{\mathcal B_{j-1}}
 Q_i(a)H'(a_i-\eta)\zeta_j\ud x
 \\
 \qquad=
 \ds\sum_{i=1}^{p}
\ds\int_{\mathcal B_{j-1}}
 \big(Q_i(a)- Q_i(1+[a-\eta]_+)\big)\ 
 \ln(1+[a_i-\eta]_+))\ \zeta_j\ud x
 \\
 \qquad\qquad+\underbrace{
  \ds\sum_{i=1}^{p}
\ds\int_{\mathcal B_{j-1}}
 Q_i(1+[a-\eta]_+)\ 
 \ln(1+[a_i-\eta]_+))\ \zeta_j\ud x
}_{\leq 0\textrm{ by h4)}}
 \\[.3cm]
\qquad 
 \leq  2p\mathscr Q\ \ds\sum_{i=1}^{p}\int_{\mathcal B_{j-1}}(1+[a_i-\eta]_+)^{q-1}\ln (1+[a_i-\eta]_+)\ud x.
 \end{array}
 \]
The last estimate is a consequence of h1) and of the elementary inequality
$$|1+[a-\eta]_+-a|\leq 1+ |[a-\eta]_+-a|\leq 1+\eta\leq 2$$ (see \cite[proof of Lemma~3.1]{GoVa} or \cite[Lemma~3]{CaVa}).
We arrive at
\[\begin{array}{l}\ds\frac{\ud}{\ud t}
\ds\sum_{i=1}^{p}\int_{\mathcal B_{j-1}}
 \zeta_jH(a_i-\eta)\ud x
 +
 4\delta_\star
 \ds\sum_{i=1}^{p}\int_{\mathcal B_{j}}\big|\nabla_x \sqrt{1+[a_i-\eta]_+}\big|^2\ud x
 \\[.3cm]
 \qquad 
 \leq C\delta^\star 2^{2j}  \ds\sum_{i=1}^{p}\int_{\mathcal B_{j-1}}
H(a_i-\eta)\ud x
\\
\qquad\qquad+
2p\mathscr Q\ \ds\sum_{i=1}^{p}\int_{\mathcal B_{j-1}}(1+[a_i-\eta]_+)^{q-1}\ln (1+[a_i-\eta]_+)\ud x.
\end{array}
\]
 We integrate this relation over $(s,t)$, with $-t_j\leq t\leq 0$ and $-t_{j-1}\leq s\leq t_j$, and next we average with respect to $s\in (-t_{j-1},-t_j)$, taking into account that
 $t_{j-1}-t_j=2^{-j-2}$.
 We obtain
 \[\begin{array}{l}
 \ds\sum_{i=1}^{p}\int_{\mathcal B_{j}}
 H(a_i-\eta)(t,x)\ud x
 +
 4\delta_\star
 \ds\sum_{i=1}^{p}\ds\int_{-t_j}^{t}\int_{\mathcal B_{j}}\big|\nabla_x \sqrt{1+[a_i-\eta]_+}\big|^2(\tau,x)\ud x\ud \tau
\\[.3cm]
 \leq 
 \ds\sum_{i=1}^{p}\int_{\mathcal B_{j-1}}
 \zeta_j H(a_i-\eta)(t,x)\ud x
 \\[.3cm]
\qquad+
 4\delta_\star
 \ds\sum_{i=1}^{p}
 \ds\frac{1}{2^{-j-2}}\ds\int_{-t_{j-1}}^{-t_j}
 \ds\int_{s}^{t}\int_{\mathcal B_{j}}\big|\nabla_x \sqrt{1+[a_i-\eta]_+}\big|^2(\tau, x)\ud x\ud\tau\ud s
\\
[.3cm]
 \leq 
\ds\sum_{i=1}^{p}\ds\frac{1}{2^{-j-2}}\ds\int_{-t_{j-1}}^{-t_j} \int_{\mathcal B_{j-1}}
 \zeta_j H(a_i-\eta)(s,x)\ud x
\ud s
\\[.3cm]
\qquad+
C\delta^\star 2^{2j}  \ds\sum_{i=1}^{p}
\ds\frac{1}{2^{-j-2}}\ds\int_{-t_{j-1}}^{-t_j} \ds\int_{s}^t
\int_{\mathcal B_{j-1}}
H(a_i-\eta)(\tau,x)\ud x\ud\tau \ud s
\\
\qquad+
2p\mathscr Q\ \ds\sum_{i=1}^{p}
\ds\frac{1}{2^{-j-2}}\ds\int_{-t_{j-1}}^{-t_j} \ds\int_{s}^t
\int_{\mathcal B_{j-1}}(1+[a_i-\eta]_+)^{q-1}\ln (1+[a_i-\eta]_+)(\tau,x)\ud x\ud \tau\ud s
\\
[.3cm]
 \leq 
 2^{j+2}
\ds\sum_{i=1}^{p}\ds\int_{-t_{j-1}}^{-t_j} \int_{\mathcal B_{j-1}}
 H(a_i-\eta)(s,x)\ud x
\ud s
\\
[.3cm]
\qquad+
C\delta^\star 2^{2j}  \ds\sum_{i=1}^{p}
 \ds\int_{-t_{j-1}}^t
\int_{\mathcal B_{j-1}}
H(a_i-\eta)(\tau,x)\ud x\ud\tau 
\\
\qquad+
2p\mathscr Q\ \ds\sum_{i=1}^{p}
\ds\int_{-t_{j-1}}^t
\int_{\mathcal B_{j-1}}(1+[a_i-\eta]_+)^{q-1}\ln (1+[a_i-\eta]_+)(\tau,x)\ud x\ud \tau.
 \end{array}\]
 We conclude by taking the supremum over $t\in(-t_j,0)$.
 %
%
%
\QED

Next, we specify the level set considered in these estimates: we use Lemma~\ref{gainofr} with 
$\eta=k_j$
and we set
$$\mathscr{U}_j= \left(\sup_{-t_j\leq t\leq 0}\sum_{i=1}^{p}\int_{{\mathcal{B}}_j} H(a_i-k_j)\ud x+
\sum_{i=1}^{p}\iint_{{\mathcal{Q}}_j}|\nabla_x\sqrt{1+[a_i-k_j]_+}|^2\ud x\ud s\right).$$


\begin{lemma}\label{unestimates} Let $2\leq q<2\frac{N+1}{N}$. 
 Then 
\begin{enumerate}[label=\roman*)]
\item For any  $ r>1$ there exists  a universal constant $c_r>0$  such that
$$\mathscr U_0\leq c_r 
 \ds\sum_{i=1}^{p} \big ( \|a_i\|^r_{L^{ r}((-1,0)\times B_1)}+ \|a_i\|^{1/2}_{L^{ r}((-1,0)\times B_1)}+
 \|a_i\|_{L^{ r}((-1,0)\times B_1)}\big ).$$

\item There exists a constant  $\Lambda>1$ such that
$$\mathscr{U}_j\leq  \Lambda^j\mathscr{U}_{j-1}^{1+N/2}$$
 for any $j\geq j_0$. Consequently, 
there exists  $\delta>0$  such that $\mathscr{U}_0\leq \delta$ implies $\lim_{j\to\infty} \mathscr{U}_j=0$.\\
 \end{enumerate}
\end{lemma}

\noindent
{\bf Proof.} Throughout the proof, we simply denote by $c$ a constant that depends 
only on the parameters of the model, and on the Lebesgue exponent, without 
paying attention to the possible changes of the value of the constant from a line to another. 

For proving i), we go back to the definition
\[\mathcal U_0=
\
\left(\ds\sup_{-\frac{1}{2}\leq t\leq 0}\ds\sum_{i=1}^{p}
\ds\int_{ \mathcal{B}_{0}}
H(a_i)\ud x
+
\ds\sum_{i=1}^{p}\ds\iint_{ \mathcal{Q}_{0}}\big|\nabla_x\sqrt{a_i+1}\big|^2\ud x\ud  \tau \right),
\]
where we remind the reader that $\mathcal B_0=B_{1/2}$ and $\mathcal Q_0=(-1/2,0)\times B_{1/2}$.
We make use of the following elementary inequalities 
\begin{align}
& H(z)\leq c \big(z (1+|\ln(z))|\big)
\label{fact2}\\
& 
|\nabla \sqrt{1+a}|\leq |\nabla\sqrt a|,
 \label{fact1}
\end{align}
which hold for any $z\geq 0$ and any  (smooth enough) function $a:\mathbb R^N\rightarrow [0,\infty)$, respectively.
We consider $\zeta_0\in C^\infty_c(\mathbb R^N)$, supported in $B_1$, such that $0\leq \zeta_0(x)\leq 1$ on $\mathbb R^N$ and 
$\zeta_0(x)=1$ on $\mathcal B_0$.
We get
\[\ds\frac{\ud}{\ud t}\ds\sum_{i=1}^{p}\ds\int_{ {B}_1} \zeta_0(x)a_i(t,x)\ud x
=
\ds\sum_{i=1}^{p}\ds\int_{{B}_1} d_i \Delta\zeta_0(x)\ a_i(t,x)\ud x
\leq \delta^\star\|\Delta \zeta_0\|_{L^\infty} \ds\sum_{i=1}^{p}\ds\int_{ {B}_1} a_i(t,x)\ud x.\]
Let $t\in(-\frac{1}{2}, 0)$ and $\tau\in (-1,t)$. We integrate over the time interval $(\tau,t)$, and then we average over $\tau\in (-1,-\frac{1}{2})$.
We are led to 
\[\begin{array}{lll}
\ds\sup_{-\frac{1}{2}\leq t\leq 0}\ds\sum_{i=1}^{p}\ds\int_{ \mathcal{B}_{0}}a_i(t,x)\ud x
&\leq&
c\ds\sum_{i=1}^{p} \ds\int_{-1}^0\ds\int_{ {B}_1}
a_i(\tau,x)\ud x\ud \tau.
\end{array}\]\\
Similarly, the localized version of the  entropy dissipation becomes
\[\begin{array}{l}
\ds\frac{\ud}{\ud t}\ds\sum_{i=1}^{p}\ds\int_{B_1} \zeta_0(x)\ a_i\ln (a_i)\ud x
+ \ds\int_{B_1} \zeta_0(x)\ds\frac{d_i |\nabla_x{a_i}|^2}{a_i}\ud x
\\[.3cm]
\qquad
=
 \ds\sum_{i=1}^{p}\ds\int_{B_1} \Delta \zeta_0 d_i\big(a_i\ln(a_i) - a_i\big)\ud x.
\\[.3cm]
\qquad
\leq \delta^\star\|\Delta \zeta_0\|_{L^\infty}
 \ds\sum_{i=1}^{p}\ds\int_{B_1} \big(a_i|\ln(a_i)| + a_i\big)\ud x.
\end{array}\]
Again we  integrate with respect to the time variable.
We shall also use the trick
\[
u|\ln(u)|=
u\ln(u)\mathbf 1_{ u\geq 1}-u\ln(u)\mathbf 1_{0\leq u< 1}\leq
u\ln(u)\mathbf 1_{ u\geq 1} +\ds\frac 2e\sqrt{u}\mathbf 1_{ 0\leq u< 1},
\]
which allows us to dominate 
\[u|\ln(u)|\leq c(u^r + \sqrt u).\]
It follows that
\[\begin{array}{l}
\ds\sup_{-\frac{1}{2}\leq t\leq 0}\ds\sum_{i=1}^{p}\ds\int_{ \mathcal{B}_{0}}a_i|\ln(a_i)|\ud x
+4\delta_\star \ds\sum_{i=1}^{p} \ds\int_{-\frac{1}{2}}^0\ds\int_{ \mathcal{B}_{0}}\big|\nabla_x\sqrt{a_i}\big|^2\ud x\ud\tau
\\
\qquad\qquad\leq
 c \ds\sum_{i=1}^{p}\left ( \ds\int_{-1}^0\ds\int_{B_1}
\left(a_i^r+\sqrt a_i + a_i\right)\ud x\ud  \tau \right)
\\[.3cm]
\qquad\qquad\leq
 c \ds\sum_{i=1}^{p}\left( \ds\int_{-1}^0\ds\int_{B_1}
|a_i|^r\ud x\ud \tau +  \left(\ds\int_{-1}^0\ds\int_{B_1} |a_i|^r\ud x\ud  \tau\right)^{1/(2r)}\mathrm{meas}(B_1)^{1-1/(2r)}\right.
\\[.3cm]
\qquad\qquad\qquad\qquad\qquad\left.
+
\left(\ds\int_{-1}^0\ds\int_{B_1} |a_i|^r\ud x\ud  \tau\right)^{1/r}\mathrm{meas}(B_1)^{1-1/r}\right)
,
\end{array}\]
by using  the H\"older inequality.
\\

We turn to the proof of ii). The estimate in Proposition~\ref{gainofr} can be recast as 
\begin{equation}\label{toto}
\begin{array}{lll}
\mathscr U_j&\leq& C\left(
2^{2j}\ds\sum_{i=1}^p\ds\iint_{\mathscr Q_{j-1}} H(a_i-k_j)(s,x)\ud x\ud s\right.
\\[.3cm]&&\qquad\qquad\left.+
\ds\sum_{i=1}^p\ds\iint_{\mathscr Q_{j-1}} (1+[a_i-k_j]_+)^{q-1}\ln (1+[a_i-k_j]_+)(s,x)\ud x\ud s\right).
\end{array}\end{equation}
Let us set
\[
\Psi(z)=
\sqrt{1+z}-1.\]
For any $\gamma\geq 1$, $\beta>0$, we can find a constant $c_{\gamma,\beta}$ such that
\[(1+z)^{\gamma}\ln(1+z)\leq c_\beta \Psi(z)^{2(\gamma+\beta)}.\]
Moreover, for $z\geq k_j\geq k_{j-1}$ we have
\[1\leq \ds\frac{z-k_{j-1}}{k_j-k_{j-1}}=2^j(z-k_{j-1}).\]
Hence, we can estimate both integrals in the right hand side of \eqref{toto}
by
an expression like
\[\begin{array}{l}
\ds\sum_{i=1}^p\ds\iint_{\mathcal Q_{j-1}} 2^{\gamma j} (1+[a_i-k_{j-1}]_+)^\gamma
\ln (1+[a_i-k_{j-1}]_+)\ud x\ud s
\\
\qquad
\leq
c_{\gamma,\beta} 2^{\gamma j}\ \ds\sum_{i=1}^p\ds\iint_{\mathcal Q_{j-1}}  \Psi([a_i-k_{j-1}]_+)^{\gamma+\beta}
\ud x\ud s.
\end{array}\]
We can play with the exponents $\gamma$ and $\beta$ for both term so that 
we obtain a common bound from above, and we arrive at 
\[
\mathscr U_j\leq c 2^{4 j}
 \ds\sum_{i=1}^p\ds\iint_{\mathcal Q_{j-1}}  \Psi([a_i-k_{j-1}]_+)^{2(N+2)/N}
\ud x\ud s.\]
This is possible as far as  $2(q-1)\leq 2\frac{N+2}{N}$ that is to say $q\leq 2\frac{N+1}{N}$.
We shall conclude by using an interpolation argument.
Indeed, on the one hand, we obviously have
\[
\ds\sup_{-t_{j-1}\leq s\leq 0}\ds\int_{\mathcal B_{j-1}} |\Psi([a_i-k_{j-1}]_+)|^2(s,x)\ud x
\leq \mathscr U_{j-1},\]
while Gagliardo-Nirenberg-Sobolev's inequality, see \cite[Theorem p.~125]{Nir}, yields
\[\begin{array}{l}
\ds\int_{-t_{j-1}}^0\left(\ds\int_{\mathcal B_{j-1}} |\Psi([a_i-k_{j-1}]_+)|^{2N/(N-2)}(s,x)\ud x\right)^{(N-2)/N}
\ud s
\\[.3cm]
\qquad\leq 
c\ds\iint_{\mathcal Q_{j-1}}|\nabla\Psi([a_i-k_{j-1}]_+)|^2(s,x)\ud x\ud s\leq c \mathscr U_{j-1}.\end{array}\]
By using the interpolation 
\[
\ds\frac{N+2}{N}=\theta\ds\frac{2N}{N-2}+2(1-\theta),
\qquad
\theta=\ds\frac{N-2}{N}\in(0,1),\]
we combine these information into
\[\begin{array}{l}
\ds\iint_{\mathcal Q_{j-1}}
 |\Psi([a_i-k_{j-1}]_+)|^{2(N+2)/N}(s,x)\ud x\ud s
 \\ [.3cm]
 \qquad\leq 
 \ds\int_{-t_{j-1}}^0 \left(\ds\int_{\mathcal B_{j-1}}  |\Psi([a_i-k_{j-1}]_+)|^{2N/(N-2)}(s,x)\ud x\right)^{\theta}
 \\[.3cm]
  \qquad \qquad \qquad \qquad \qquad\times
 \left(\ds\int_{\mathcal B_{j-1}}  |\Psi([a_i-k_{j-1}]_+)|^{2}(s,x)\ud x\right)^{1-\theta}\ud s
\\[.3cm]
\qquad\leq 
\mathscr U_{j-1}^{1-\theta}\ds\int_{-t_{j-1}}^0
 \left(\ds\int_{\mathcal B_{j-1}}  |\Psi([a_i-k_{j-1}]_+)|^{2N/(N-2)}(s,x)\ud x\right)^{(N-2)/N}\ud s
 \leq c\mathscr U_{j-1}^{1+2/N}.
\end{array}\]
We conclude by applying Lemma~\ref{2bete}.
\QED

Once we know that $\lim_{j\to\infty} \mathscr U_j=0$
we deduce that 
\[
\ds\lim_{j\to \infty} \ds\frac{1}{t_j}\ds\sum_{i=1}^p\ds\iint_{-\mathcal Q_j} H(a_i-k_j)\ud x\ud t
=0\geq 
4\ds\sum_{i=1}^p\ds\int_{-1/4}^0\ds\int_{B(0,1/2)} H(a_i-1)\ud x\ud t.\]
It implies that 
$0\leq a_i(t,x)\leq 1$ holds for a.~e.~$(t,x)\in (-1/4,0)\times B(0,1/4)$.

\section{Proof of Proposition~\ref{P:CaVa}}
\label{app:CaVa}

It is worth giving some hints for the proof of  Proposition~\ref{P:CaVa}, which 
is fully detailed in \cite[Proposition~11 \& Corollary~12]{CaVa}.
Again, the proof heavily relies on duality arguments.
The main step consists in showing that 
\begin{equation}\label{tata}
\|\Phi(t,\cdot)\|_{L^\infty(\mathbb R^N)}\leq \|\Phi(0,\cdot)\|_{L^\infty(\mathbb R^N)}.\end{equation}
Indeed, we remind the reader that 
$\Phi(t,x)$ is determined by the convolution formula (for $N>2$)
\[
\Phi(t,x)=-C_N\ds\int_{\mathbb R^N} \frac{M(t,y)}{|x-y|^{N-2}}\ud y
\]
where $$C_N=\ds\frac{1}{(N-2)\sigma_N},$$
with $\sigma_N=\frac{2\pi^{N/2}}{\Gamma(N/2)}$ the measure of the unit sphere of $\mathbb R^N$.
Thus,  given $R>0$, we simply split 
\[
\Phi(0,x)=-C_N\ds\int_{|x-y| \leq R}\frac{M(0,y)}{|x-y|^{N-2}}\ud y
-
C_N\ds\int_{|x-y| > R} \frac{M(0,y)}{|x-y|^{N-2}}\ud y
\]
which yields
\[
|\Phi(0,x)|\leq
C_N\|M(0,\cdot)\|_{L^\infty(\mathbb R^N)}
\ds\frac{\sigma_N R^2}{2}
+\ds\frac{C_N}{R^{N-2}}\|M(0,\cdot)\|_{L^1(\mathbb R^N)}
.\]
Optimizing with respect to $R$, we get
\[
|\Phi(0,x)|\leq K_N \|M(0,\cdot)\|_{L^\infty(\mathbb R^N)}^{1-2/N}\|M(0,\cdot)\|_{L^1(\mathbb R^N)}^{2/N},
\]
where $K_N>0$ depends only on the space dimension $N\geq 3$.
\\

In order to justify \eqref{tata}, we need  to introduce a mollified diffusion coefficient. Indeed, as the $a_i$'s are smooth on $[0,T_{\mathrm{max}})\times \mathbb R^N$,
$M$ is smooth too; thus $(t,x)\mapsto d(t,x)$ is a smooth function, except possibly at the points where  $M(t,x)$ vanishes.
Given $\mu>0$, we denote $d_\mu(t,x)$ a smooth function verifying
\[d_\mu(t,x)=d(t,x)\quad \textrm{when $M(t,x)\geq \mu$},\qquad
0< \delta_\star\leq d_\mu(t,x)\leq \delta^\star.\]
The proof of \eqref{tata} splits into two steps.
 
Let $0<T<\infty$, Let $\zeta\in C^\infty_c(\mathbb R^N)$ and consider the solution of the \emph{final} equation
\begin{equation}\label{tata2}
\begin{array}{l}
\partial_t \varphi + d_\mu \Delta \varphi=0
,
\qquad
\varphi(T,x)=\zeta(x),
\end{array}\end{equation}
together with the initial value problem
\[
\partial_t \rho - \Delta (d_\mu \rho)=0
,
\qquad
\rho(0,x)=\rho^0(x).
\]
We assume that
\[\|\zeta\|_{L^\infty(\mathbb R^N)}\leq 1.\]
The maximum principle, see for instance \cite[Theorem~8, Chapter~7]{Evans},
implies 
\[\ds\sup_{0\leq t\leq T}\|\varphi(t,\cdot)\|_{L^\infty(\mathbb R^N)}\leq\|\zeta\|_{L^\infty(\mathbb R^N)}\leq 1
.\]
We have
\[
\ds\frac{\ud}{\ud t}\ds\int_{\mathbb R^N}\rho (t,x)\varphi(t,x)\ud x=0.\]
It follows that
\[\left|
\ds\int_{\mathbb R^N}\rho(T,x)\zeta(x)\ud x\right|
=
\left|
\ds\int_{\mathbb R^N}\rho^0(x)\varphi(0,x)\ud x\right|
\leq \|\rho^0\|_{L^1(\mathbb R^N)}.\]
By virtue of the Hahn-Banach theorem, we conclude that 
\[\begin{array}{lll}
\|\rho(T,\cdot)\|_{L^1(\mathbb R^N)}
&=&\sup\left\{\left|
\ds\int_{\mathbb R^N}\rho (T,x)\zeta(x)\ud x\right|,\
\zeta\in C^\infty_c(\mathbb R^N),\,
\|\zeta\|_{L^\infty(\mathbb R^N)}\leq 1\right\}
\\[.3cm]
&
\leq& \|\rho^0\|_{L^1(\mathbb R^N)}.\end{array}\]

Next, we shall apply a similar reasoning in order to make the norm 
$
\|\Delta \zeta\|_{L^1(\mathbb R^N)}$ appear.
For $0<T<\infty$ and $\varphi$ solution of \eqref{tata2}, let us set $$\rho(t,x)=\Delta\varphi(T-t,x)$$
which satisfies
$$
\partial_t \rho-\Delta (d_\mu \rho)=0,\qquad \rho(0,x)=\Delta \zeta(x)\in L^1(\mathbb R^N).$$
The previous step thus tells us that
\[
\|\rho(T,\cdot)\|_{L^1(\mathbb R^N)} =
\|\Delta \varphi(0,\cdot)\|_{L^1(\mathbb R^N)}\leq \|\rho(0,\cdot)\|_{L^1(\mathbb R^N)}=
\|\Delta \zeta\|_{L^1(\mathbb R^N)}.\]
Going back to the equation for the total mass, we get
\[\ds\frac{\ud}{\ud t}\ds\int_{\mathbb R^N}
M\varphi (t,x)\ud x=\ds\int_{\mathbb R^N}
M(d-d_\mu)\Delta\varphi (t,x)\ud x.\]
Let $0<T<T_{\mathrm{max}}$. Integrating over $(0,T)$  yields 
\[\begin{array}{l}
\left|
\ds\int_{\mathbb R^N}
M\varphi (T,x)\ud x
\right|
= 
\left|
\ds\int_{\mathbb R^N}
\Delta \Phi\varphi (T,x)\ud x
\right|
=
\left|\ds\int_{\mathbb R^N}
 \Phi(T,x)\Delta\zeta (x)\ud x
\right|
\\
[.3cm]
\qquad
=
\left|
\ds\int_{\mathbb R^N}
M\varphi (0,x)\ud x+
\ds\int_0^T\ds\int_{\mathbb R^N}
M(d-d_\mu)\Delta\varphi (t,x)\ud x\ud t
\right|
\\[.3cm]
\qquad
=
\left|
\ds\int_{\mathbb R^N}
\Delta \Phi\varphi (0,x)\ud x
+
\ds\int_0^T\ds\int_{\mathbb R^N}
M(d-d_\mu)\Delta\varphi (t,x)\ud x\ud t
\right|
\\[.3cm]\qquad
\leq 
\left|
\ds\int_{\mathbb R^N}
 \Phi \Delta \varphi (0,x)\ud x
\right|+
\left|+
\ds\int_0^T\ds\int_{\mathbb R^N}
M(d-d_\mu)\Delta\varphi (t,x)\ud x\ud t
\right|
\\
[.3cm]
\qquad\leq 
\|\Phi(0,\cdot)\|_{L^\infty(\mathbb R^N)}\|\Delta \varphi(0,\cdot)\|_{L^1(\mathbb R^N)}
+2T\delta^\star \mu  \|\Delta\varphi\|_{L^\infty(0,T;L^1(\mathbb R^N)}
\\[.3cm]
\qquad \qquad\qquad
\textrm{since 
$ |d-d_\mu|M=|d-d_\mu| M\mathbf 1_{M\leq \mu}\leq 2\delta^\star \mu$}
\\
[.3cm]
\qquad
\leq 
\big(\|\Phi(0,\cdot)\|_{L^\infty(\mathbb R^N)}
+2T\delta^\star \mu \big)
\|\Delta \zeta\|_{L^1(\mathbb R^N)}.\end{array}\]
This relation holds for any $\mu>0$ and $\zeta\in C^\infty_c(\mathbb R^N)$.
Therefore, we can conclude that \eqref{tata} holds, which ends the proof.
\QED

\section*{Acknowledgements}
T.~G. acknowledges warm welcome of  the Math.~Department of UT-Austin
where a couple of visits have made possible 
progress on this question. A.~V. is partially supported  by the NSF grant DMS-1614918.

\bibliography{CaGoVa}
\bibliographystyle{plain}

\end{document}